\documentclass[plain]{ipart}

\RequirePackage[OT1]{fontenc}
\RequirePackage{amsthm,amsmath}
\RequirePackage[numbers,square]{natbib}
\RequirePackage[colorlinks]{hyperref}

\usepackage{fullpage}
\usepackage{graphicx}
\usepackage{amssymb}
\usepackage{epstopdf}
\usepackage{caption, subcaption}
\usepackage{float}
\usepackage{bbm}
\newcommand{\D}{\mathcal{D}}
\DeclareMathOperator{\sd}{sd}

\RequirePackage[pagewise]{lineno} %line number package

\newcommand*{\QED}{\hfill\ensuremath{\square}\\}

\startlocaldefs
\theoremstyle{plain}
\newtheorem{thm}{Theorem}[section]
\newtheorem{cor}[thm]{Corollary}
\newtheorem{claim}[thm]{Claim}
\newtheorem{rem}[thm]{Remark}
\newtheorem{lem}[thm]{Lemma}
\newtheorem{conj}[thm]{Conjecture}

\endlocaldefs

\begin{document}

\setpagewiselinenumbers

\begin{frontmatter}
\title{Toward \.{Z}ak's conjecture on graph packing\protect}
%\runtitle{}

\begin{aug}
\author{\fnms{Ervin} \snm{Gy\H{o}ri}
\thanksref{t1}
\ead[label=e1]{ervin@renyi.hu}},
\address{Alfr\'{e}d R\'{e}nyi Institute of Mathematics \\ Budapest, Hungary \\ and\\ Department of Mathematics, Central European University \\ Budapest, Hungary
\printead{e1}}

\author{\fnms{Alexandr} \snm{Kostochka}\thanksref{t2}\ead[label=e2]{kostochk@math.uiuc.edu}}
\address{Department of Mathematics \\ University of Illinois \\ Urbana, IL 61801, USA\\and\\ Sobolev Institute of Mathematics\\ Novosibirsk, Russia
\printead{e2}}

\author{\fnms{Andrew} \snm{McConvey}
\ead[label=e3]{mcconve2@illinois.edu}}
\address{Department of Mathematics \\ University of Illinois \\ Urbana, IL 61801\\ USA
\printead{e3}}

\and

\author{\fnms{Derrek} \snm{Yager}\ead[label=e4]{yager2@illinois.edu}}
\address{Department of Mathematics \\ University of Illinois \\ Urbana, IL 61801\\ USA
\printead{e4}}

\thankstext{t1}{Research of this author is supported in part by OTKA Grants 78439 and 101536}
\thankstext{t2}{Research of this author is supported in part by NSF grant  DMS-1266016 and  by grants 12-01-00631 and 12-01-00448 of the Russian Foundation for Basic Research.}
\runauthor{E. Gy\H{o}ri et al.}

\end{aug}

\begin{abstract}
Two graphs $G_{1} = (V_{1}, E_{1})$ and $G_{2} = (V_{2}, E_{2})$, each of order $n$, \emph{pack} if there exists a bijection $f$ from $V_{1}$ onto $V_{2}$ such that  $uv \in E_{1}$ implies $f(u)f(v) \notin E_{2}$.  In 2014, \.{Z}ak proved that if $\Delta (G_{1}), \Delta (G_{2}) \leq n-2$ and $|E_{1}| + |E_{2}| + \max \{ \Delta (G_{1}), \Delta (G_{2}) \} \leq 3n - 96n^{3/4} - 65$, then $G_{1}$ and $G_{2}$ pack.  In the same paper, he conjectured that if $\Delta (G_{1}), \Delta (G_{2}) \leq n-2$, then $|E_{1}| + |E_{2}| + \max \{ \Delta (G_{1}), \Delta (G_{2}) \} \leq 3n - 7$ is sufficient for $G_{1}$ and $G_{2}$ to pack.  We prove that, up to an additive constant, \.{Z}ak's conjecture is correct.  Namely, there is a constant $C$ such that if $\Delta(G_1),\Delta(G_2) \leq n-2$ and $|E_{1}| + |E_{2}| + \max \{ \Delta(G_{1}), \Delta(G_{2}) \} \leq 3n - C$, then $G_{1}$ and $G_{2}$ pack.  In order to facilitate induction, we prove a stronger result on list packing.
\end{abstract}

\dedicated{Dedicated to Adrian Bondy on the occasion of his $70^{\text{th}}$ birthday.}



\end{frontmatter}

%======================================================================
%Section:  Introduction
%======================================================================
% 
%
\section{Introduction}\label{sec:intro}
%
%
%======================================================================
%======================================================================
Extremal problems on graph packing have been actively  studied since the seventies. 
Recall that  two $n$-vertex  graphs  are said to \emph{pack} if there is an edge-disjoint placement of the graphs onto the same set of vertices.  
More technically, a \emph{packing} of graphs $G_1$ and $G_2$ is a bijection $f:V_1\rightarrow V_2$ such that for all $u,v\in V_1,$ either $uv\notin E_1$ or $f(u)f(v)\notin E_2$. 
In 1978, Bollob\'{a}s and Eldridge~\cite{B-E} and  Sauer and Spencer~\cite{S-S} proved several important results on graph packing. 
In particular, Sauer and Spencer~\cite{S-S} showed that two $n$-vertex graphs  pack if the product of their maximum degrees is less than $n/2$.

\begin{thm}[\cite{S-S}]\label{S-S product}
Let $G_1$ and $G_2$ be two $n$-vertex  graphs. If $2\Delta(G_1) \Delta(G_2) < n$, then $G_1$ and $G_2$ pack.
\end{thm}

For $n=2k$ with $k$ odd, if $G_1=K_{k,k}$ and $G_2$ is a perfect matching $M_k$, then $G_1$ and $G_2$ do not pack; so the bound is sharp.  Sauer and Spencer~\cite{S-S} and Bollob\'{a}s and Eldridge~\cite{B-E} independently proved sufficient conditions for packing two graphs with given average degrees.

\begin{thm}\label{S-S}
Let $G_1$ and $G_2$ be two $n$-vertex  graphs. If $|E(G_{1})| + |E(G_{2})| \leq \frac{3}{2}n - 2$ then $G_1$ and $G_2$ pack.
\end{thm}

Moreover, Bollob\'{a}s and Eldridge proved that Theorem~\ref{S-S} can be significantly strengthened when we additionally assume that $\Delta(G_{1}), \Delta(G_{2}) < n-1$.

\begin{thm}[\cite{B-E}]\label{B-E} Let $G_1$ and $G_2$ be two $n$-vertex  graphs. 
If $\Delta(G_1),\Delta(G_2) \leq n-2, |E(G_1)| + |E(G_2)| \leq 2n-3,$ and $\{G_1,G_2\}$ is not one of the following pairs: $\{2K_2, K_1\cup K_3\}, \{\overline{K_2} \cup K_3, K_2 \cup K_3\}, \{3K_2, \overline{K_2} \cup K_4\}, \{\overline{K_3} \cup K_3, 2K_3\},$ $\{2K_2 \cup K_3,\overline{K_3}\cup K_4\},\{\overline{K_4} \cup K_4, K_2 \cup 2K_3\}, \{\overline{K_5} \cup K_4, 3K_3\}$, then  $G_1$ and $G_2$ pack.
\end{thm}

This theorem is also sharp, which we see by observing that graphs $G_1=K_{1,n-2} \cup K_{1}$ and $G_2=C_{n}$ do not pack.  
Recently, \.{Z}ak showed that with stronger restrictions on maximum degrees of $G_1$ and $G_2$ we can weaken restrictions on their sizes.  Namely, he proved the following.

\begin{thm}[\cite{Z}]\label{thm:zak 5/2}
Let $G_{1}$ and $G_{2}$ be two graphs of order $n \geq 10^{10}$.  If $|E(G_{1})| + |E(G_{2})| + \max \{ \Delta(G_{1}), \Delta(G_{2}) \} < \frac{5}{2}n - 2$, then $G_{1}$ and $G_{2}$ pack.
\end{thm}

 \.{Z}ak showed that this result can also be strengthened by forbidding the star on $n$ vertices.

\begin{thm}[\cite{Z}]\label{zak}
Let $G_1$ and $G_2$ be $n$-vertex graphs with $\Delta (G_{1}), \Delta (G_{2}) \leq n-2$.   If $|E(G_{1})| + |E(G_{2})| + \max \{ \Delta(G_{1}), \Delta(G_2) \} \leq 3n - 96n^{3/4}-65$, then $G_{1}$ and $G_{2}$ pack.
 \end{thm}
 
This theorem is asymptotically sharp since $K_{1,n-2} \cup K_{1}$ and $C_{n}$ do not pack. In the same paper \.{Z}ak poses the following conjecture.

\begin{figure}[h!]
\begin{center}
\includegraphics[width=.6\textwidth]{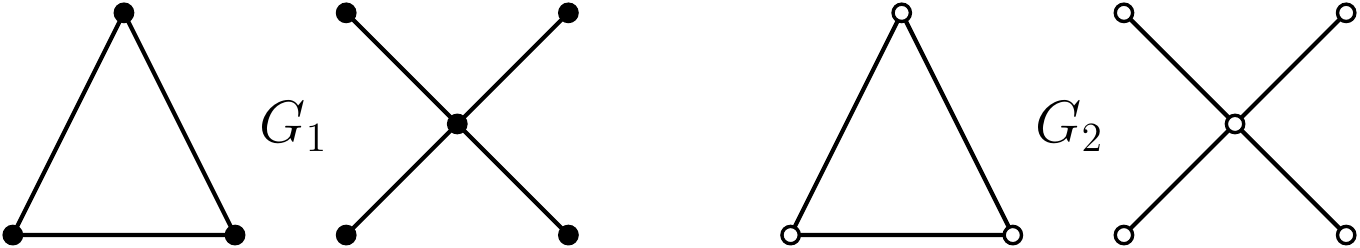}
\end{center}
\caption{Sharpness example for Conjecture~\ref{conj:zak}.  In this example $n=8$ and $|E(G_{1})| + |E(G_{2})| + \max{\{ \Delta(G_{1}), \Delta(G_2)\} } = 3n-6$ but the graphs do not pack.}
\label{fig:zaksharp}%
\end{figure}

\begin{conj}[\cite{Z}]\label{conj:zak}
Let $G_1$ and $G_2$ be $n$-vertex graphs with $\Delta (G_{1}), \Delta(G_{2}) \leq n-2$.   If $|E(G_{1})| + |E(G_{2})| + \max{\{ \Delta(G_{1}), \Delta(G_2)\} } \leq 3n - 7$, then $G_{1}$ and $G_{2}$ pack.
\end{conj}

\.{Z}ak also provides the following example to show that, if true, the conjecture is best possible.   Let $n \geq 8$ and let $G_{1}$ and $G_{2}$ each be isomorphic to $K_{3} + K_{1, n-4}$, a disjoint union of a triangle and a star (Figure~\ref{fig:zaksharp}).  Then, $\Delta(G_1)  = \Delta(G_2) = n-4$ and $|E(G_{1})| + |E(G_{2})| + \max{\{ \Delta(G_{1}), \Delta(G_2)\} } = (n-1) + (n-1) + (n-4) = 3n-6$.  A simple check shows that $G_{1}$ and $G_{2}$ do not pack.

\begin{figure}[h]
\begin{center}
\includegraphics[width=.6 \textwidth]{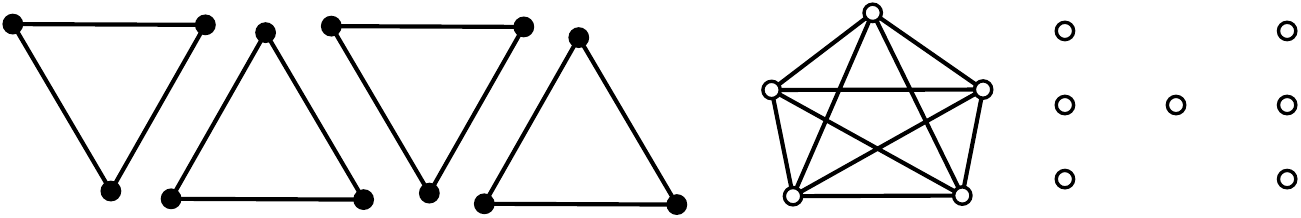}
\end{center}
\caption{\.{Z}ak's Conjecture is false for small values of $n$.}
\label{fig:ZakConj}%
\end{figure}

However, for some small values on $n$,  Conjecture~\ref{conj:zak} fails. 
 For example, consider $G_1 = 4 K_3$ and $G_2=K_5\cup\overline{K}_7$ (Figure~\ref{fig:ZakConj}). 
In any attempted packing, we are forced to send at least two vertices from the same component in $G_1$ to the clique in $G_2$, so the graphs do not pack. 
 In this example, $|E(G_1)| + |E(G_{2})|+ \max{\{ \Delta (G_1), \Delta(G_2) \}} = 12 + 10 + 4 = 26 = 3n - 10$. 
 We were unable to find large counterexamples, so the conjecture may hold with a finite set of exceptions.  
Further, the main result of this paper shows that, up to the choice of the additive constant, Conjecture~\ref{conj:zak} is true.

\begin{thm}\label{thm:main}
Let $C = 11(195^2) = 418,275$. Let $G_1$ and $G_2$ be $n$-vertex  graphs  with $\Delta (G_{1}), \Delta(G_{2}) \leq n-2$. If $|E(G_{1})| + |E(G_{2})| + \max{\{ \Delta(G_{1}), \Delta(G_2)\} } \leq 3n - C$, then $G_{1}$ and $G_{2}$ pack.
\end{thm}

Our constant $C$ is not optimal and we can somewhat decrease it by a more detailed case analysis in our proofs.  However, since $3n - 96n^{3/4}-65 \leq 0$ for $n \leq 10^{6}$, Theorem~\ref{thm:main} improves the previous best known result even for small values of $n$.  Further, Theorems~\ref{thm:main} and~\ref{S-S} together imply that Theorem~\ref{thm:zak 5/2} holds when $n$ is at least $2C - 2 \approx 10^{6}$.  To see this notice that if $\Delta(G_{1}) = n - 1$ or $\Delta(G_{2}) = n-1$, then $|E(G_{1})| + |E(G_{2})| \leq \frac{3}{2}n -1$ and Theorem~\ref{S-S} applies.  Alternatively, when $n \geq 2C-2$, $\frac{5}{2}n - 2 \leq 3n - C$ and Theorem~\ref{thm:main} applies.

Our proof of Theorem~\ref{thm:main} uses the concept of list packing introduced in \cite{ListB-E}.  A {\em graph triple} $\mathbf{G} = (G_{1}, G_{2}, G_{3})$ consists of two disjoint $n$-vertex graphs $G_1=(V_1,E_1)$ and $G_2=(V_2,E_2)$ and a bipartite graph $G_{3} = (V_{1} \cup V_{2}, E_{3})$
with partite sets $V_1$ and $V_2$.  A \emph{list packing} of $\mathbf{G}$ is a packing of $G_{1}$ and $G_{2}$ such that $uf(u)\notin E_3$ for any $u \in V_1$. Essentially, a list packing is a packing of $G_{1}$ and $G_{2}$ with an additional set of restrictions on the bijection $f$.

We prove the following list version of Theorem~\ref{thm:main}.

\begin{thm}\label{thm:list main} Let $C=11(195^2)$.  Let $n \geq 2$ and $\mathbf{G}=(G_1,G_2,G_3)$ be a graph triple with $|V_1|=|V_2|=n$, $\Delta (G_1),\Delta (G_2) \leq n-2$, 
and $\Delta (G_3) \leq n-1$.  
If $|E_{1}| + |E_{2}| + |E_{3}| + \max \{ \Delta(G_{1}), \Delta (G_{2}) \} + \Delta(G_{3}) \leq 3n - C$, then $\mathbf{G}$ packs.
\end{thm}

\begin{figure}[h!]
\begin{subfigure}{.32\textwidth}
  \centering
  \includegraphics[width=.5\linewidth]{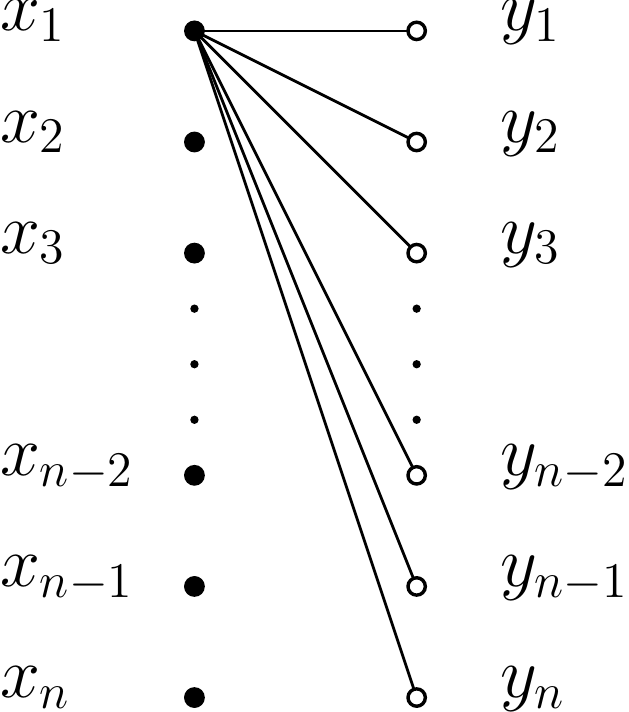}
  \label{fig:yellow star}
  \caption{}
\end{subfigure}%
\begin{subfigure}{.32\textwidth}
  \centering
  \includegraphics[width=.5\linewidth]{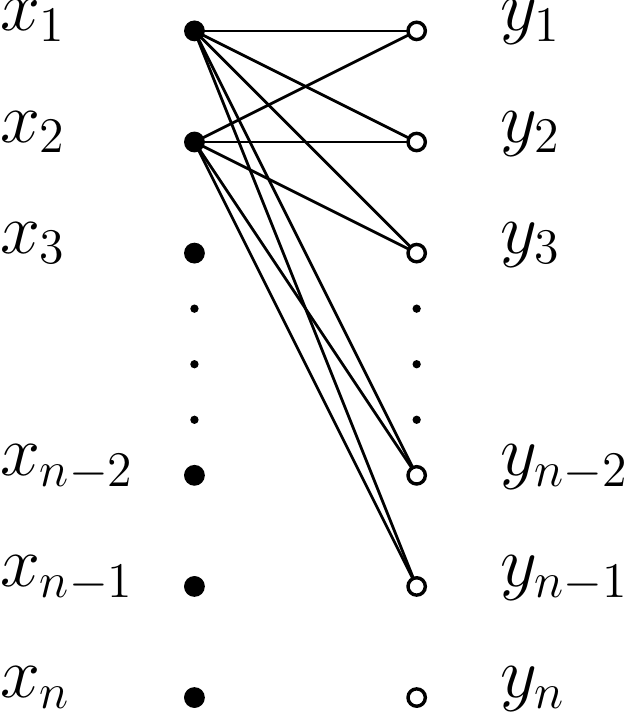}
  \label{fig:double yellow star}
  \caption{}
\end{subfigure}
\begin{subfigure}{.32\textwidth}
  \centering
  \includegraphics[width=.5\linewidth]{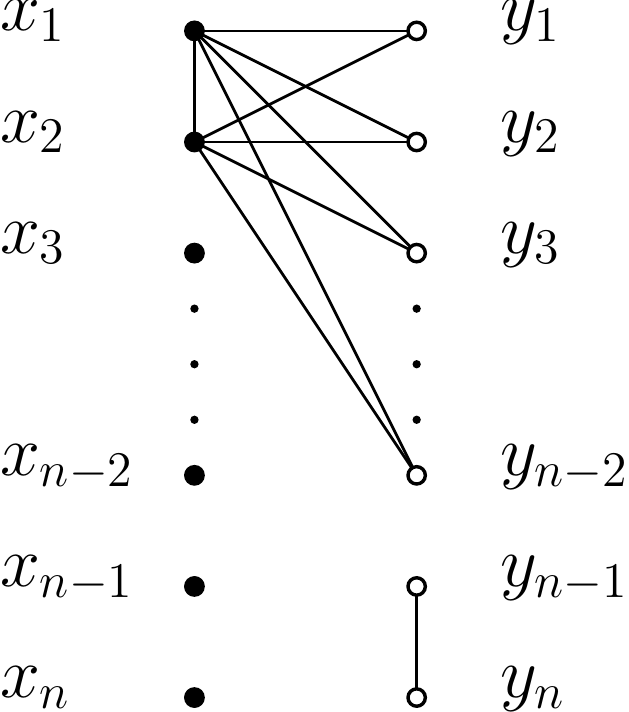}
  \label{fig:two edges}
  \caption{}
\end{subfigure}

\begin{subfigure}{.32\textwidth}
  \centering
  \includegraphics[width=.5\linewidth]{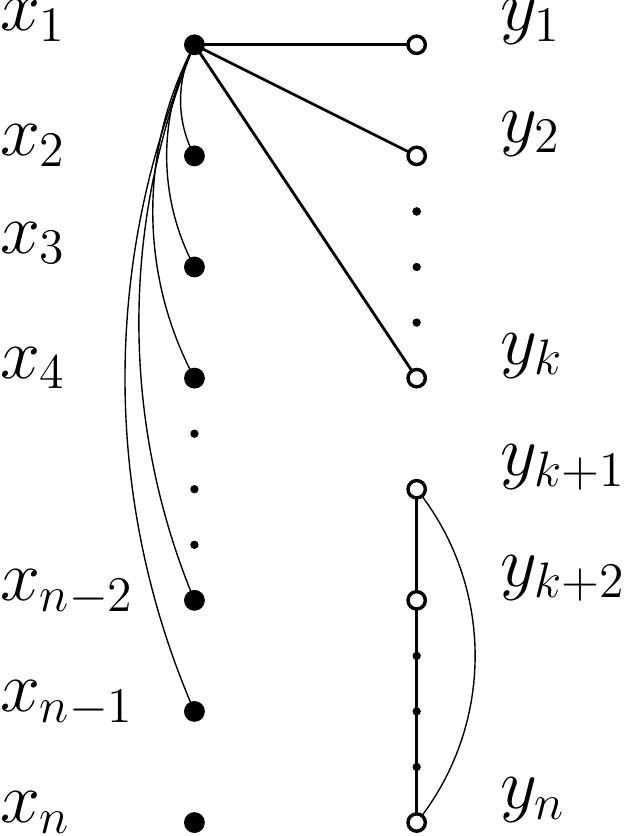}
  \label{fig:partition}
  \caption{}
\end{subfigure}
\begin{subfigure}{.32\textwidth}
  \centering
  \includegraphics[width=.5\linewidth]{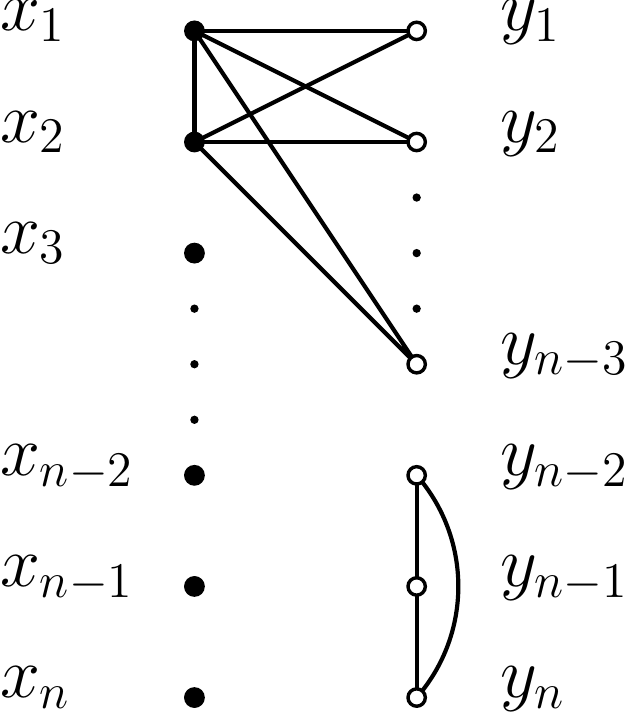}
  \label{fig:shortcycle}
  \caption{}
\end{subfigure}
\begin{subfigure}{.32\textwidth}
  \centering
  \includegraphics[width=.5\linewidth]{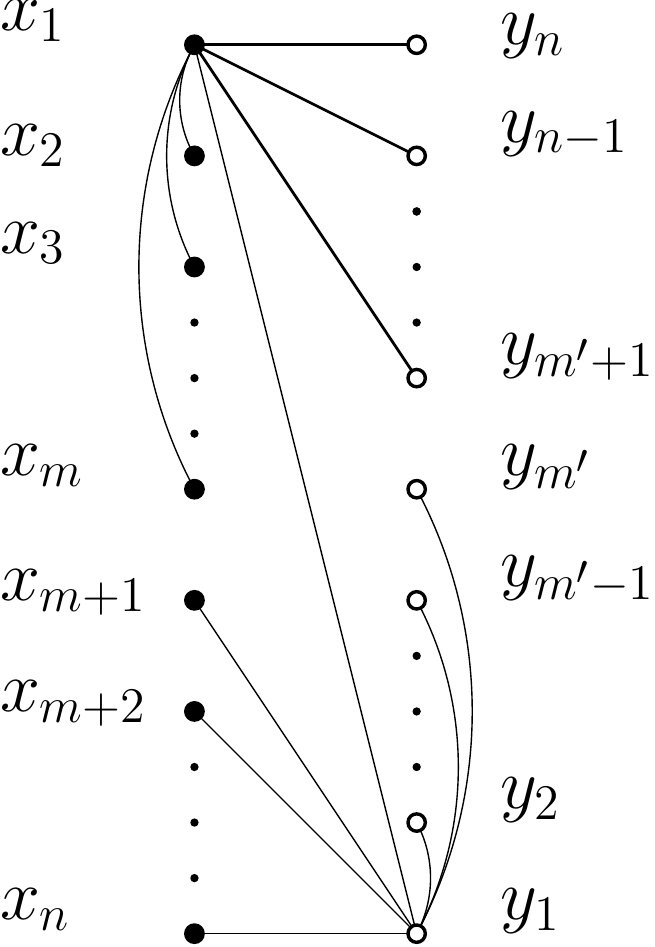}
  \label{fig:cycle}
  \caption{}
\end{subfigure}
\caption{Sharpness examples for Theorem~\ref{thm:list main}}
\label{fig:sharpness}
\end{figure}

Note that Theorem~\ref{thm:main} is the special case of Theorem~\ref{thm:list main} in which $G_{3}$ has no edges. 
The pair shown in Figure~\ref{fig:ZakConj} shows that, up to an additive constant, the theorem is sharp.  
Moreover, there are other infinite families of examples showing that, up to an additive constant, the theorem is sharp  when $E_{3}$ has linear
in $n$ number of edges.  Several of these examples are shown in Figure~\ref{fig:sharpness}.  
The body of this paper  contains a proof of the slightly stronger Theorem~\ref{thm:detailed main}. 
 This theorem is more technical than Theorem~\ref{thm:list main} and we refer the reader to Section~\ref{sec:prelim} 
for the statement of the theorem and an explanation of the used notation.

The paper is organized as follows.  In Section~\ref{sec:prelim}, we state definitions, some useful preliminary results, and the main technical result,
Theorem~\ref{thm:detailed main}.  The proof of Theorem~\ref{thm:detailed main} will be by contradiction. In Section~\ref{sec:min counterexample} we prove several lemmas regarding the degree requirements of a minimal counterexample $\mathbf{G}=(G_1,G_2,G_3)$.  We then use these properties in Section \ref{sec:1-neighbors} to show that a minimal counterexample has at most one vertex with at least two neighbors of degree $1$. 
Next, in Section \ref{sec:sponsors}, we introduce the notion of  supersponsors and show that each of $G_1$ and $G_2$ contains at least two supersponsors.  Finally, in Section~\ref{sec:end}, we arrive at a contradiction by using the structure of a minimal  counterexample to construct a packing.

%======================================================================
%Section:  Preliminaries
%======================================================================
% 
%
\section{The setup}\label{sec:prelim}
%
%
%======================================================================
%======================================================================

A graph triple $\mathbf{G} = (G_{1}, G_{2}, G_{3})$ of order $n$ consists of a pair of $n$-vertex graphs $G_{1} = (V_1, E_1)$ and $G_{2} = (V_{1}, E_{2})$ together with a bipartite graph $G_{3} = (V_{1} \cup V_{2}, E_{3})$.  Let $V(\mathbf{G}):= V_{1} \cup V_{2}$ be the vertex set of the graph triple, $E(\mathbf{G})=E_1 \cup E_2 \cup E_3$ be the edge set of the graph triple, and $e(\mathbf{G})=|E(\mathbf{G})|$. We omit $\mathbf{G}$ when it is clear.  The triple $\mathbf{G}$ \emph{packs} if there is a bijection $f:V_{1} \rightarrow V_{2}$ such that $v f(v) \notin E_{3}$ for any $v \in V_{1}$ and $uv \in E_{1}$ implies $f(u)f(v) \notin E_{2}$.  An edge in $E_{1} \cup E_{2}$ is a \emph{white} edge, while an edge in $E_{3}$ is a \emph{yellow} edge.

For $v \in V_{i}$ ($i = 1,2$), the \emph{white neighborhood} of $v$, denoted $N_{i}(v) \subseteq V_{i}$, is the set of neighbors of $v$ in $G_{i}$ and $d_{i}(v) = |N_{i} (v)|$.  For convenience, when $w \in V_{3-i}$, we say that $N_{i}(w) = \emptyset$ (and hence $d_{i}(w) = 0$).  The \emph{yellow neighborhood} of $v \in V_{i}$, denoted $N_{3}(v) \subseteq V_{3-i}$ is the set of neighbors of $v$ in $G_{3}$ and $d_{3} (v) = |N_{3}(v)|$.  Vertices in the white (respectively, yellow) neighborhood of $v$ are called \emph{white neighbors} (respectively, \emph{yellow neighbors}).  For $v \in V_{i}$, the \emph{neighborhood} of $v$, denoted $N(v)$ is the disjoint union $N_{i}(v) + N_{3}(v)$ and the \emph{degree} of $v$ is $d_{i}(v) + d_{3}(v)$ and is denoted $d(v)$. Also, we use $N[v]$ to denote the \emph{closed neightborhood} of $v$, i.e. $N[v]=N(v)\cup \{v\}$. For disjoint vertex sets $X$ and $Y$ in a graph triple, $\|X,Y\|$  denotes the number of edges connecting  $X$ and $Y$.
For brevity, if $X=\{x\}$ and $Y=\{y\}$, then we will write $\|x,y\|$ instead of $\|\{x\},\{y\}\|$.

 When considering a specific graph triple $\mathbf{G}$, we will let $e_{i} = |E_{i}|$ and define $\Delta_{i} = \max_{v \in V} d_{i}(v)$ for $i = 1,2,3$.
In \cite{ListB-E}, the authors proved extensions of Theorem~\ref{S-S product} and Theorem~\ref{B-E} to list packing.  The following two theorems will be used throughout this paper.

\begin{thm}[\cite{ListB-E}]\label{List S-S product}
 Let $\mathbf{G} = (G_{1}, G_{2}, G_{3})$ be a graph triple with $|V_{1}| = |V_{2}| = n$.  If $\Delta_{1} \Delta_{2} + \Delta_{3} \leq n/2,$ then $\mathbf{G}$ does not pack if and only if $\Delta_{3} = 0$ and one of $G_1$ or $G_2$ is a perfect matching and the other is $K_{\frac{n}{2},\frac{n}{2}}$ with $\frac{n}{2}$ odd or contains $K_{\frac{n}{2}+1}.$ Consequently, if $\Delta_{1} \Delta_{2} + \Delta_{3} < n/2,$ then $\mathbf{G}$ packs.\end{thm}

\begin{thm}[\cite{ListB-E}]\label{List B-E} Let $\mathbf{G} = (G_{1}, G_{2}, G_{3})$ be a graph triple with $|V_{1}| = |V_{2}| = n$.  If $\Delta_{1}, \Delta_{2} \leq n-2$, $\Delta_{3} \leq n-1$, $|E_{1}| + |E_{2}| + |E_{3}| \leq 2n-3$ and the pair $(G_1,G_2)$ is none of the 7  pairs in  Theorem~\ref{B-E}, then $\mathbf{G}$ packs.
\end{thm}

\noindent For a graph triple $\mathbf{G} = (G_{1}, G_{2}, G_{3})$, let $\Delta_{3|i} = \max_{v \in V_{i}} d_{3} (v)$, $D_{i} = \max  \{ \Delta_{i}, \Delta_{3|i} \}$, and  \[\D = \max{ \{ \Delta_{1} + \max\{ \Delta_{3|2} - 4, 0 \}, \Delta_{2} +\max\{ \Delta_{3|1} - 4, 0\} \} }.\]
Instead of Theorem~\ref{thm:list main}, it is more convenient to prove the following.  

\begin{thm}\label{thm:detailed main} Let $C:=11(195^2) + 4$. Let $n \geq 2$ and $\mathbf{G}=(G_1,G_2,G_3)$ be a graph triple of order $n$.  If
\begin{equation}
\label{e0}
\mbox{  $\Delta_1,\Delta_2 \leq n-2$, $\Delta_3 \leq n-1$}
\end{equation}
and
\begin{equation}
\label{e1}
\quad F(\mathbf{G}): =  e_1+ e_2 + e_3 + \D \leq 3n - C,
\end{equation}
then $\mathbf{G}$  packs.
\end{thm}

Note that Theorem~\ref{thm:detailed main} implies Theorem~\ref{thm:list main} since $\Delta_{3} \geq \Delta_{3|1}, \Delta_{3|2}$ and $F(\mathbf{G}) + 4 \leq e_{1} + e_{2} + e_{3} + \max \{ \Delta_{1}, \Delta_{2} \} + \Delta_{3}$.  In proving this theorem, we will often consider two graph triples, $\mathbf{G}$ and $\mathbf{G}'$ and will compare $F(\mathbf{G})$ and $F(\mathbf{G}')$.  Define $\partial (\mathbf{G},\mathbf{G}') = F(\mathbf{G}) - F(\mathbf{G}')$.  The rest of the paper will be a proof of Theorem~\ref{thm:detailed main}.

%======================================================================
%Section: Minimal Counterexample
%======================================================================
% 
%
\section{Maximum and Minimum Degrees in a Minimal Counterexample}\label{sec:min counterexample}
%
%
%======================================================================
%======================================================================
Fix $C:=11(195^2) + 4$ and let $\mathbf{G} = (G_{1}, G_{2}, G_{3})$ be a graph triple of the smallest order $n$ such that $\mathbf{G}$ satisfies~\eqref{e0} and~\eqref{e1} but $\mathbf{G}$ does not pack.  By Theorem~\ref{List B-E} and~\eqref{e1},
 \begin{equation}
\label{e2}
\D \leq n+2-C.
\end{equation}
This yields $n\geq C-2$. Moreover, since $n\geq C-2$, Theorem~\ref{List S-S product} implies $\D\geq 2$, and thus, by~\eqref{e2}, $n\geq C$.

\begin{lem}\label{3l}
Every vertex of $\mathbf{G}$ has a white neighbor.
\end{lem}

\noindent\textbf{Proof:} Suppose $v\in V$ has no white neighbor. Without loss of generality, let $v \in V_1$.

\textbf{Case 1:} \emph{The vertex $v$ is isolated in $\mathbf{G}$.} If any $w \in V_2$ has degree at least $3$ in $\mathbf{G}$ then  taking $\mathbf{G}'=(G_1-v,G_2-w,G_3-v-w)$ and $n'=n-1$ 
gives $\partial (\mathbf{G},\mathbf{G}')\geq3$ and thus $F(\mathbf{G}')\leq 3n'-C$. Also by~\eqref{e2}, for $i=1,2$, 
\[\Delta_i'\leq\Delta_i\leq\D+4\leq n+6-C\leq (n-1)-2=n'-2.\] So by the minimality of $\mathbf{G}$,
 the new triple $\mathbf{G}'$ packs. Then  this packing extends to a packing of $\mathbf{G}$ by sending $v$ to $w$, contradicting the choice
of  $\mathbf{G}$.  So suppose the  degree of each $w\in V_2$ is at most $2$.  By Theorem~\ref{List S-S product}, there is a vertex $v' \in V_{1}$ with $d(v')>n/6$. By~\eqref{e0},  there  is a non-neighbor $w$ of $v'$ in $V_2$.  If $w$ has a white neighbor, say $y\in V_2$, then let $\mathbf{G}''=(G_1-v-v',G_2-w-y,G_3-v-v'-w-y)$ with $n''=n-2$; otherwise, let $\mathbf{G}''=(G_1-v',G_2-w,G_3-v'-w)$ with $n''=n-1$.
Then  $\partial (\mathbf{G},\mathbf{G}'')>d(v')=n/6>6$ and so $F(\mathbf{G}'')\leq 3n''-C$ which by~\eqref{e2} implies $\Delta_i''\leq n+6-C\leq n''-2$ for $i=1,2$, .
Thus again by the minimality of $\mathbf{G}$, the triple
 $\mathbf{G}''$ packs.  Then, we extend this packing of $\mathbf{G}''$ to a packing of $\mathbf{G}$ by  sending $v'$ to $w$ (and $v$ to $y$ if $y$ exists), again contradicting the choice of $\mathbf{G}$.
 
The last subcase of Case 1 is that  $d_2(w)=2$ for every non-neighbor $w$ of $v'$ in $V_2$. In particular, $e_2+e_3\geq e_2+ d_{3}(v') \geq n$.  
So, if $X=V_1-N[v']-v$, then by \eqref{e1}
\[\sum_{x\in X}d_1(x)\leq 2e_1-2d_1(v')\leq 2\left[ 3n - C - \D - (e_2 + d_{3}(v')) - d_1(v') \right].\]
Since $d_1(v')+|X|=n-2$, $e_{3} \geq d_{3}(v')$, and $\D \geq \Delta_1\geq d_1(v')$, we get
\[ \sum_{x\in X}d_1(x)\leq 2\left( 3n - C - 2d_1(v') - n \right)\leq 2(2|X|+4-C) < 4|X| -8.\]

So, there are nonadjacent $x_1,x_2\in X\subset V_1$ with $d_1(x_1),d_1(x_2)\leq 3$.

\begin{figure}[H]
\begin{center}
\includegraphics[width=.2 \textwidth]{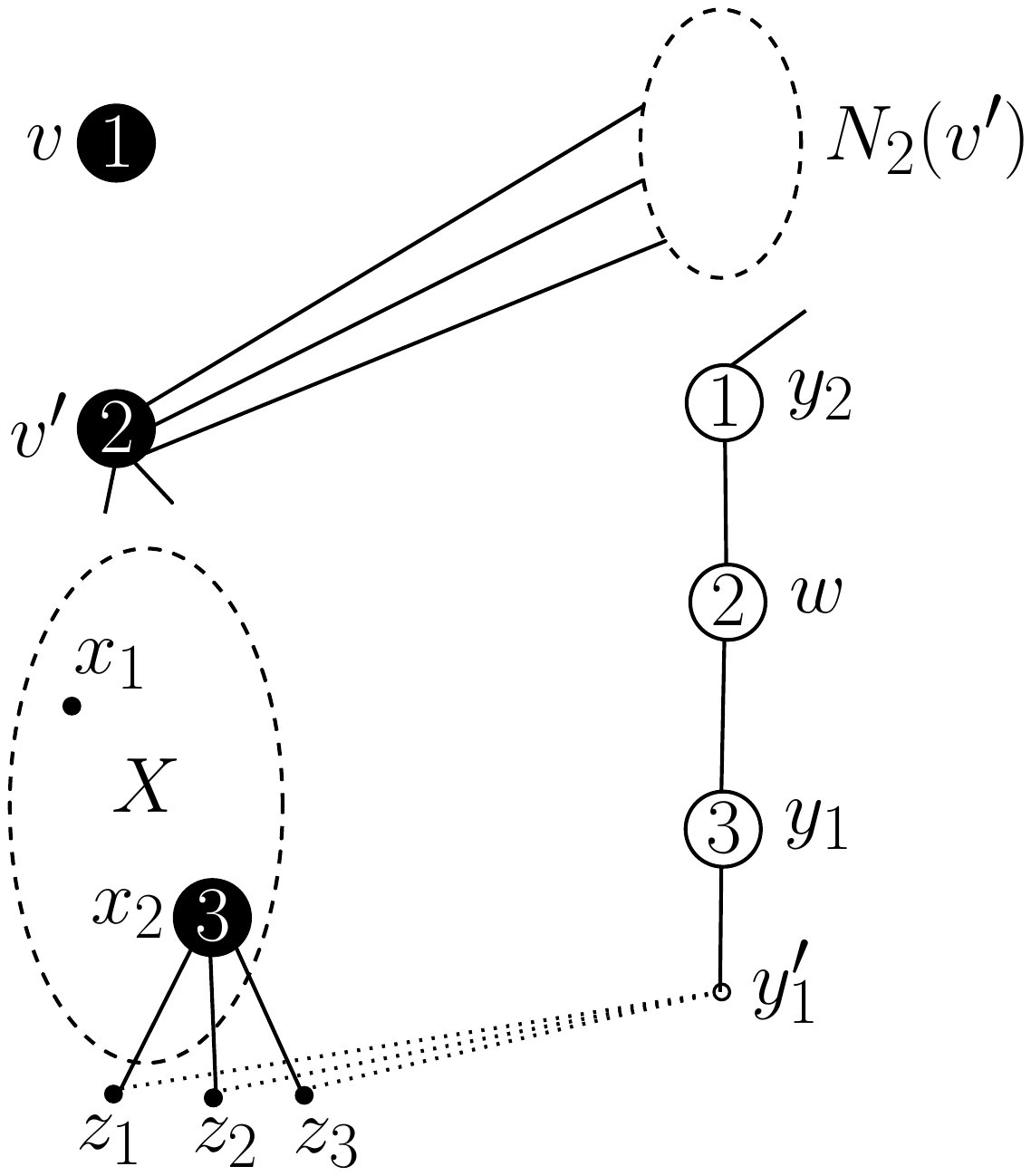}
\end{center}
\caption{Packing used at the end of Case 1}
\label{fig:Lemma11}%
\end{figure}

Let $w$ be a non-neighbor of $v'$ in $V_2$ and let $y_1$ and $y_2$ be the white neighbors of $w$.
 Since $y_1w\in E_2$ and $d(y_1)\leq 2$, we may assume $y_1x_2\notin E_3$.
 Choose $z_1,z_2,z_3\in V_1$ so that $N_1(x_2)\subset \{z_1,z_2,z_3\}$.  Let $y'_1$ be the white neighbor of $y_1$ distinct from $w$, if exists.
  Then we place $v'$ on $w$, $v$ on $y_2$, $x_2$ on $y_1$, and add yellow edges from $y'_1$ to $N_1(x_2)$ (Figure~\ref{fig:Lemma11}).
 Since this decreases $e_1+e_2+e_3$ by at least $n/6+2\geq C/6+2\geq 12$  and increases $\D$ by at most $3$,  we are left with a graph triple $\mathbf{G}'$ of order at least $n-3$ 
 and $F(\mathbf{G}') \leq 3(n-3) - C$.  Also by~\eqref{e2}, both inequalities in~\eqref{e0} hold.
 So by the minimality of $\mathbf{G}$, there is a packing of $\mathbf{G}'$, and this packing extends of a packing of $\mathbf{G}$.

\textbf{Case 2:} \emph{The vertex $v \in V_1$ is incident to yellow edges.}  Let $A := N_{3} (v)$. By the case, $|A| \geq 1$. Since $V_2-A\neq\emptyset$ by \eqref{e2}, there is some $w \in V_2-A$. Since Case 1 does not hold, $d(w)\geq 1$.
 If $d(v)+d(w)\geq 3$, then we can construct a packing by sending $v$ to $w$ and creating a new graph triple $\mathbf{G}'$ by removing these two vertices.  
 In creating $\mathbf{G}'$, we have removed $3$ edges, and observe that by~\eqref{e2}, the inequalities in~\eqref{e0} holds for $\mathbf{G}'$.  So $\mathbf{G}'$ packs by the minimality of $\mathbf{G}$,
  and this packing extends to a packing of the original triple, a contradiction.
 Thus, $d(v)=1$ (say $A=\{w'\}$) and $d(w)=1$ for each $w\in V_2-w'$.
 
 Let $Y=V_2-N[w']$. Since $d_2(w')\leq\Delta_{2} \leq \D \leq n + 2 - C$, we have $|Y|\geq C-3$.  If $d(w')=1$, then by switching the roles of 
 $v$ and $w'$, we conclude that $d(v')=1$  for each $v'\in V_1-v$; so $\mathbf{G}$ packs by
 Theorem~\ref{List S-S product}. Hence, $d(w')\geq 2$.  There are two cases.
 
\textbf{Case 2.1:} \emph{$G_2[Y]$ has no edges.} Since the white neighbors of $w'$ cannot have other neighbors, every $y \in Y$ has no white neighbors.  
If also  every vertex in $V_1$ has degree $1$, then by~\eqref{e2}, 
\[ e_{1} + e_{2} + e_{3} = \frac{(2n-1) +d(w')}{2} \leq n - \frac{1}{2} + \D+4 \leq n - \frac{1}{2} + (n + 6 - C) < 2n - 3.
\]
In this case, $\mathbf{G}$ packs by Theorem~\ref{List B-E}, a contradiction. So we conclude that there is a vertex $x \in V_{1}$ of degree at least $2$.  

Next, assume that two vertices $y_{1}, y_{2} \in Y$ have distinct neighbors in $V_1$.  Then we may assume that $x$ is not adjacent to one of these vertices, say $y_{1}$, 
and let $\mathbf{G}'=(G_1-x,G_2-y_1,G_3-x-y_1)$ and $n'=n-1$. 
 Since $\partial (\mathbf{G},\mathbf{G}')\geq 3$  and~\eqref{e0} holds for $\mathbf{G}'$ by~\eqref{e2}, 
 $\mathbf{G}'$ packs by the minimality of $\mathbf{G}$, and this packing extends to a packing of $\mathbf{G}$ by placing $x$ on $y_1$.
 
 Hence, each vertex in $Y$ is adjacent to the same vertex $x' \in V_{1}$.  This implies $\D \geq d_{2}(w') + d_{3} (x') - 4 \geq n - 5$, a contradiction to ~\eqref{e2}.

\textbf{Case 2.2:} \emph{There is an edge $y_1y_2 \in E(G_2[Y])$.} Then
 \begin{equation}
\label{dec24}
\mbox{for every non-adjacent $x_1,x_2\in V_1$, $d(x_1)+d(x_2)\leq 4$, }
\end{equation}
since otherwise we could send $x_1$ to $y_1$ and $x_2$ to $y_2$ and consider $\mathbf{G}''=(G_1-x_1-x_2,G_2-y_1-y_2,G_3-x_1-x_2-y_1-y_2)$.
We have $\partial (\mathbf{G},\mathbf{G}'')\geq 6$  and~\eqref{e0} holds for $\mathbf{G}''$ by~\eqref{e2}, so
 $\mathbf{G}''$ packs by the minimality of $\mathbf{G}$, and this packing extends to a packing of $\mathbf{G}$. 

Since none of $x\in V_1-v$ is  adjacent to $v$, by~\eqref{dec24}, $d(x) \leq 3$ for every $x\in V_1$,   
In particular, this yields $\Delta_{1} \leq 3$,  $\Delta_{2}= \max\{1,d_2(w')\} \leq 1+d_2(w')$, and $\Delta_{3}\leq \max\{3,d_3(w')\}\leq 3+d_3(w')$. 
Then, \[\Delta_1\Delta_2+\Delta_3\leq 3(d_2(w')+1)+(3+d_3(w')) \leq 3(d(w')+2).\] Since $\mathbf{G}$ does not pack, Theorem~\ref{List S-S product} implies that $\Delta_1\Delta_2+\Delta_3\geq n/2$, so  $d(w') \geq \frac{n}{6}-2$.

By \eqref{e2}, $n+2-C\geq \D \geq d_3(w')-4$, so there are at least $C-6$ non-neighbors of $w'$ in $V_1$.  
By~\eqref{dec24},   at most $4$ vertices in $V_{1}$ have degree $3$. Thus there exists a non-neighbor $x_0$ of $w'$ such that $d(x_0) \leq 2$ and the degrees of the white neighbors of $x_0$, which could be neighbors of $w'$, as well, also do not exceed $2$. If $N_1(x_0)=\emptyset,$ then send $x_0$ to $w'$.  If $N_{1}(x_0) = \{z_1\}$, then send $x_0$ to $w'$, $z_1$ to $y_1$ and $v$ to $y_2$.
 If $N_{1}(x_0) = \{z_1,z_2\}$ and $z_1z_2\notin E_1$, then send $x_0$ to $w'$, $z_1$ to $y_1$ and $z_2$ to $y_2$.
 Finally, if $N_{1}(x_0)= \{z_1,z_2\}$ and $z_1z_2\in E_1$, then by the choice of $x_0,z_1,z_2$,
 these 3 vertices induce a component in $\mathbf{G}$; so we can send $x_0$ to $w'$, $z_1$ to $y_1$ and $z_2$ to any $y_0\in Y-y_2$.  In all cases, we have deleted at least $\frac{n}{6} - 2$ edges.  Since by~\eqref{e2},~\eqref{e0} also will hold in all cases, we can pack the resulting graph triple, and then extend this to a packing of $\mathbf{G}$, a contradiction. \QED

\begin{lem}\label{4l}
If a vertex in $V_1$ has degree $1$, then no vertex in $V_2$ has degree $1$.
\end{lem}

\noindent\textbf{Proof:} Suppose $v\in V_1, w\in V_2$ and $d(v)=d(w)=1$. Then by Lemma~\ref{3l}, the edges incident to $v$ and $w$ are white.
 Let $vv'\in E_1$ and $ww'\in E_2$. Let $A_1=N_1(v')-v$, $A_2=N_3(v')=N(v')\cap V_2$,  $B_1=N_3(w')=N(w')\cap V_1$, $B_2=N_2(w')-w$.
  Let $x_0$ (respectively, $y_0$) be a vertex of maximum degree among the vertices in $V_1-v-v'$ (respectively, in $V_2-w-w'$).
   
   We   obtain graph triple $\mathbf{G}'=(G_1',G_2',G_3')$ by first placing $v'$ on $w$, $v$ on $y_0$, deleting the matched pairs, and then 
   adding yellow edges from $w'$ to the vertices in $A_1 \setminus B_{1}$. 
    If $\mathbf{G}'$ packs, then together with our placement of $v'$ on $w$ and $v$ on $y_0$ we will have a packing of $\mathbf{G}$.  If it does not pack, then by the minimality of $\mathbf{G}$, 
    either~\eqref{e0} or~\eqref{e1} does not hold for $\mathbf{G}'$.  Since $\Delta_1, \Delta_{2} \leq \D \leq n - C + 2$ and the white degrees of vertices did not increase, if~\eqref{e0} is violated in $\mathbf{G}'$, then by~\eqref{e2}, $\mathbf{G}'$ has a vertex $u$ with $d_3'(u) = n-2$.  Since $\Delta_3 = \max\{ \Delta_{3|1}, \Delta_{3|2} \} \leq \D + 4$,~\eqref{e2} implies that $u=w'$.  However, $n-2 \leq d_{3}'(w') \leq d_{1}(v') + d_{3}(w') \leq \Delta_{1} + \Delta_{3|2} \leq \D + 4$,  a contradiction to~\eqref{e2}. Thus~\eqref{e1} must be violated in $\mathbf{G}'$:
\begin{equation}
\label{12la}
F(\mathbf{G}') = e(G_1') + e(G_2') + e(G_3')+\D' \geq 3(n-2)-C+1.
\end{equation} 
 Symmetrically, we   obtain graph triple $\mathbf{G}''=(G_1'',G_2'',G_3'')$ by first placing   $v$ on $w'$ and $x_{0}$ on $w$, deleting the matched pairs, and then
  adding yellow edges from $v'$ to the vertices in $B_{2} \setminus A_{2}$.  Similarly to~\eqref{12la}, we derive
\begin{equation}
\label{12lb}
F(\mathbf{G}'') = e(G_1'') + e(G_2'') + e(G_3'') +\D'' \geq 3(n-2)-C+1.
\end{equation} 

The proof also will require the following claim.

\begin{claim}\label{claim:DegreeBounds}
If there exist constants $a, b$ such that $d(x_{0}) \leq a$, $d(y_{0}) \leq b$, and $C-3 \geq \max \{ 2a(b+2), 2(a+2)b \}$, then $\mathbf{G}$ packs.
\end{claim}

\noindent\textbf{Proof of Claim:}  By symmetry, we will assume that $a \geq b$ so that $C-3 \geq 2a(b+2)$.   We will construct a packing of $\mathbf{G}$ that maps $v$ to $y_{0}$, $v'$ to $w$.  Observe that since $|A_{1}| + |B_{1}| \leq (\Delta_{1} - 1) + \Delta_{3|2} \leq \D + 3 \leq n - C + 5$, we may choose a vertex $x \in V_{1} - N_{1}[v'] - N_{3} [w']$ that we may map to $w'$.  In order to preserve the packing property, we must ensure that white neighbors of $x$ are not mapped to white neighbors of $w'$.  Again, by~\eqref{e2}, we see that there are at least $C-3$ vertices of $V_{2} - N_{2}[w']$.  Since $y_{0}$ has maximum degree among all vertices in $V_{2} - w'$, the average degree of the vertices in this set is at most $b$.  By Turan's Theorem, we may find an independent set of vertices in $V_{2} - N_{2}[w']$ of size at least $(C-3) / (b + 1) \geq 2a$.

Now, let $\{x_{1}, \ldots, x_{a'}\} = N_{1} (x)$ be the white neighborhood of $x$ and notice that $a' = d_{1}(x) \leq d(x_{0}) \leq a$.  Since $x_{0}$ was maximal, $d_{3}(x_{i}) \leq a - 1$, for each $i = 1, \ldots, a'$.  Thus, we may successively map each $x_{i}$ on a non-neighbor $y_{i}$ chosen from the independent set in $V_{2} - N_{2} [w']$.  After each such mapping, we add yellow edges between the white neighbors of $x_{i}$ and the white neighbors of $y_{i}$.  This yields a new graph triple $\mathbf{G}^*$ of order $n - a'  - 3$.  In this new triple, we see that $\Delta_{1}^{*} \leq a, \Delta_{2}^{*} \leq b$ and, due to the added yellow edges, $\Delta_{3}^{*} \leq a + b - 2$.  However, this gives 
\[2\Delta_{1}^{*}\Delta_{2}^{*} + 2\Delta_{3}^{*} \leq 2ab + 2(a + b - 2) \leq 2ab + 4a \leq C \leq n- a' - 3.\]
By Theorem~\ref{List S-S product}, $\mathbf{G}^{*}$ packs and this packing extends to a packing of $\mathbf{G}$. \QED

Along with this Claim, we will use~\eqref{12la} and~\eqref{12lb}  to prove the lemma.   Observe that to obtain $\mathbf{G}'$, we deleted $|A_{1}| + |A_{2}| + 1$ edges adjacent to $v'$, one edge adjacent to $w$, $d(y_{0})$ edges adjacent to $y_{0}$ (though we may have double counted the edge $v'y_{0}$), and added $|A_{1} \setminus B_{1}|$ new yellow edges adjacent to $w'$.  Thus, by \eqref{12la} and similarly by \eqref{12lb}, 
\begin{equation}\label{G'}
5 \geq \partial (\mathbf{G},\mathbf{G}') \geq |A_{1} \cap B_{1}| + |A_{2}| + d(y_{0}) + 1 + \D - \D'.
\end{equation}

\begin{equation}\label{G''}
5 \geq \partial (\mathbf{G},\mathbf{G}'') \geq |A_{2} \cap B_{2}| + |B_{1}| + d(x_{0}) + 1 + \D - \D''.
\end{equation}
If $\D - \D' \geq -1$ and $\D - \D'' \geq -1$, then $d(x_{0}), d(y_{0}) \leq 5$ and we are done by Claim~\ref{claim:DegreeBounds}.  So by symmetry, we may assume  that $\D - \D'' \leq -2$.  In particular, since the only vertex in $\mathbf{G}''$ that has increased its degree by more than 1 is $v'$, we  have  $\D'' = \Delta_{2}'' + d_{3}''(v') - 4$.  There are two cases.

\textbf{Case 1:} \emph{$\D - \D' \leq -2$.}  In creating $\mathbf{G}'$, the only vertex that has increased its degree by at least $2$ is $w'$, so $\D' = \Delta_{1}' + d_{3}'(w') - 4$.  Observing that $d_{3}'(w') = |A_{1} \cup B_{1}|$ and plugging this in for $\D'$ and $\D''$, we can sum together~\eqref{G'} and~\eqref{G''} to get
\begin{equation}\label{dec241}
10 \geq 2|A_{1} \cap B_{1}| + 2|A_{2} \cap B_{2}| + d(y_{0})+ d(x_{0}) + 2\D - \Delta_{1}' -\Delta_{2}'' - |A_{1}| - |B_{2}| + 10. 
\end{equation}
Since $\D \geq \Delta_{1}, \Delta_{2}$,  we have $\D \geq |A_{1}| + 1$ and  $\D \geq |B_{2}| + 1$.  Furthermore, since $x_{0}$ was a maximum degree vertex in $V_{1} - v'$, we have $d(x_{0}) \geq \Delta_{1}'$. Similarly, $d(y_{0}) \geq \Delta_{2}''$.  Inserting these inequalities into~\eqref{dec241}, we get
\[ 10 \geq   2|A_{1} \cap B_{1}| + 2|A_{2} \cap B_{2}| + 12.\]
This is a contradiction, so the case is proved.

\textbf{Case 2:} \emph{ $\D - \D' \geq -1$.}  We see from~\eqref{G'}  that $5 \geq |A_{1} \cap B_{1}| + |A_{2}| + d(y_{0})$.  Note also, that since $w'$ is a vertex in $\mathbf{G}'$, $|B_{2}| \leq d'_{2}(w') + 1 \leq \D' - \Delta_{3|1}'+ 5 \leq \D - \Delta_{3|1}' + 6$.  Next, observe that $d_{3}''(v') \leq |A_{2} \cup B_{2}|$, so we have \[\D'' \leq \Delta_{2}'' + |B_{2}| + |A_{2} \setminus B_{2}|  - 4 \leq \Delta_{2}'' + \D+ |A_{2} \setminus B_{2}|  - \Delta_{3|1}' +2.\]  We now substitute these inequalities into~\eqref{G''},
\begin{align*}
5 &\geq |A_{2} \cap B_{2}| + |B_{1}| + d(x_{0})+ 1 + \D -  \Delta_{2}'' - \D - |A_{2} \setminus B_{2}|  + \Delta_{3|1}' - 2 \\
&\geq 2|A_{2} \cap B_{2}| + |B_{1}| + d(x_{0}) -  \Delta_{2}'' - |A_{2}| + \Delta_{3|1}' -1.
\end{align*}
However, $y_{0}$ is a vertex in $\mathbf{G}''$, so $\Delta_{2}'' \leq d(y_{0}) + 1$.  In particular, 
\[ d(y_{0}) + |A_{2}| + 7 \geq 2|A_{2} \cap B_{2}| + |B_{1}| + d(x_{0}) + \Delta_{3|1}' .\]
Finally, recall that $\D - \D' \geq -1$ implies by~\eqref{G'} that $5 \geq |A_{1} \cap B_{1}| + |A_{2}| + d(y_{0})$.  This  gives that $d(y_{0}) \leq 5$, and when combined with the last inequality, that $d(x_{0}) \leq 12$.   Since $C>1000$, by  Claim~\ref{claim:DegreeBounds},  $\mathbf{G}$ packs, a contradiction.
\QED

From now on, by Lemma~\ref{4l}, we will assume that
\begin{equation}\label{dec243}
 \mbox{$d(w)\geq 2$ for every $w\in V_2$.}
\end{equation}

\begin{lem}\label{5l}
$D_1,D_2\geq 3$.
\end{lem}

\noindent\textbf{Proof:} Suppose $D_2 \leq 2$, the case where $D_1 \leq 2$ follows similarly.  The white components of $G_{2}$ are paths and cycles.  By Theorem~\ref{List S-S product}, $D_1\geq n/6$.  Also, by \eqref{e1}, \[\sum_{v\in V_1} d(v)+2\D \leq 6n-2C-\sum_{w\in V_2}d(w)  < 5n-2C.\]  
Let $v' \in V_1$ have maximum degree in $V_1$, so that $d(v')\geq n/6$. Since $\D\geq D_1-4$, this implies 
\begin{equation}\label{e7}
\sum_{v \in V_1 - \{v'\}} d(v)\leq 5n-2C-d(v')-2\D\leq 5n-2C-n/6-2(n/6-4) < 9n/2 -2C+8.
\end{equation}

Consider a vertex $w_0 \in V_{2} - N_{3}(v')$. There are two cases.

\textbf{Case 1:} \emph{The white component containing $w_0$ is not a triangle.}  In this case, $w_0$ has at most two white neighbors, $w_1, w_2 \in V_{2}$.  (Notice $w_2$ may not exist).  Since $D_{2} \leq 2$, there are at most $4$ vertices of $V_{1} - N_{1} [v']$ adjacent to $N_2(w_0)$.  By \eqref{e7}, there are at most $60$ vertices of degree at least $n/12-6$ in $V_{1} - N[v']$. So, there are at least two vertices in $V_1-N[v']$ that have degree less than $n/12-6$ and are not adjacent to $N(w_0)$, call them $v_1, v_2$. 
 We will map $v'$ to $w_0$, $v_1$ to $w_1$, and (if $w_2$ exists) $v_2$ to $w_2$.  Create a new triple $\mathbf{G}' = (G_1', G_2', G_3')$ by deleting these matched pairs and adding new yellow edges from $(N_{1}(v_{1}) - v_{2})$ to $(N_{2}(w_{1}) - w')$ and $(N_{1}(v_{2}) - v_{1})$ to $(N_{2}(w_{2}) - w')$. 
  Since $\mathbf{G}'$ has order at least $n-3$ and $\D \leq n - C +2$, we see that \eqref{e0} holds for $\mathbf{G}'$.  Notice that $w_{i}$ has at most one white neighbor other than $w'$, so we have added at most $d_{1}(v_{1}) + d_{1}(v_{2})$ new yellow edges.
Thus, $\mathbf{G}'$ has at most $e_1 + e_2 + e_3 - d(v') - d(v_1) - d(v_2) + d_1 (v_1) + d_1(v_2)$ edges and $\D' \leq \D + d_1(v_1) + d_1 (v_2)$.  
Finally, since $d(v_{i}) \geq d_{1}(v_{i})$, we have 
\begin{equation}
\label{e3}
e_1' + e_2' + e_3'+\D' \leq e_1 + e_2 + e_3 + \D - ( d (v') - d_1(v_1) - d_1(v_2) ).
\end{equation} 
\vspace{.1in}

If $e_1' + e_2' + e_3' + \D' \leq 3(n-3) - C$, then $\mathbf{G}'$ packs by the minimality of $\mathbf{G}$ and this packing extends to a packing of $\mathbf{G}$.  
But we have chosen $v_1$ and $v_2$ so that $d(v_1), d(v_2) < n/12 - 6$.  Since $d(v') \geq n/6$, we have $d (v') - d_1(v_1) - d_1(v_2) \geq 9$ and, by \eqref{e3}, $\mathbf{G}'$ packs and this extends to a packing of $\mathbf{G}$, a 
contradiction.

\textbf{Case 2:} \emph{The white component containing $w_0$ is a triangle.}  Let $w_0 w_1 w_2$ be a triangle in $G_2$ and let $d = d_{1} (v').$ Note that $d \leq \D < n -  C + 2$.  As before, there are at most $4$ vertices in $V_1 - N_1 [v']$ adjacent to $\{ w_1, w_2 \}$.  Let $X = V_{1} - N_1[v'] - N_3(\{w_1, w_2 \})$ and notice that $|X| \geq n - d - 5 \geq C - 7$.  If there are nonadjacent vertices  $x_1, x_2 \in X$, then we can match $v'$ to $w_0$, $x_1$ to $w_1,$ and $x_2$ to $w_2$.  Since $d(v') \geq n/6$, removing these vertices leaves a smaller graph triple which we can pack by the minimality of $\mathbf{G}$.  This packing extends to a packing of $\mathbf{G}$, a contradiction.

On the other hand, if all vertices of $X$ are adjacent to each other, then there are at least $\binom{|X|}{2} \geq 2 |X|$ edges in $G_{1}[X]$.  Since $v'$ has $d$ white neighbors, we see that $e_{1} + \D \geq 2|X| + 2d \geq 2n - 10$.  Finally, $e_{2} + e_{3} \geq \frac{1}{2} \sum_{w \in V_{2}} d(w) \geq n$.  So, $e_{1} + e_{2} + e_{3} + \D \geq 3n - 10$, a contradiction.
\QED

\begin{lem}\label{6l}
$\D+\sum_{v\in V_1} d(v)\geq 2n-12$.
\end{lem}

\noindent\textbf{Proof:} The sum of degrees of vertices in a component $M$ of $G_1$ containing a cycle is at least $2|V(M)|$. Thus if
$\sum_{v\in V_1} d(v)< 2n-12$, then $G_1$ has at least six tree-components, each adjacent to at most one yellow edge. Let $H$ be a smallest such component and $vw$ be the yellow edge incident to $V(H)$, if it exists. Then $s:=|V(H)|\leq n/6$. Let $w_1\in V_2$ with the maximum white degree and begin by finding a permissible vertex $v_{1}$ to send to $w_{1}$.  If $vw$ does not exist, then choose $v_1$ to be any vertex in $V(H)$.  If $vw$  exists and $w_1\neq w$, then choose $v_1=v$. Finally, if  $vw$  exists and $w_1= w$, then choose $v_1$ to be any vertex in $V(H)-v$.  Consider $H$ as a rooted tree with root $v_1$, so that each $x\in V(H)-v_1$ has a unique parent in $H$.  Order the vertices of $H$: $v_1,\ldots,v_s$ in the Breadth-First order.  We now will consecutively place all vertices of $H$ on vertices in $V_2$. We start by placing $v_1$ on $w_1$. Then for every $i=2,\ldots, s$, if possible, we place $v_i$ on a vertex  $w_i\in V_2$ not adjacent to the image $w_{i'}$ of any $v_{i'}$ with $i'<i$, and if not possible, then just on any non-occupied non-neighbor of the image $w_j$ of its parent $v_j$.

First, we show that we always can choose a vertex to place each $v_i$. Indeed, otherwise for some $2\leq i\leq s$, we cannot place $v_i$ and let's call its parent $v_j$.  Then, each vertex of $V_2$ either is adjacent to $w_j$ or is occupied by one of $v_1,\ldots,v_{i-1}$. If $j=1$, then
because $H$ is a tree obtained via Breadth-First search, $i\leq d_1(v_1)+1$. Thus in this case,
$d_2(w_1)+d_1(v_1)\geq n-1$ and since $v_{1} \in H$, $d_{2}(w_{1}) \geq \frac{3}{4}n$.  But then 
\[\D+\sum_{v\in V_1} d(v) \geq d_{2}(w_{1}) + \left(d_{1}(v_{1}) + \sum_{v \in V_1 - v_{1}} d(v)\right) \geq 2n - 2,\] contradicting our assumption. Otherwise,
 the host, say $w_j\neq w_1$, of the parent $v_j$ of $v_i$ has at least $n-i+1$ neighbors in $V_2$. Then by the choice of $w_1$, also $\D \geq d_{2}(w_1) \geq n-i+1$. Thus the total number of edges incident to $w_1$ and $w_j$ is at least $d(w_1)+d(w_j)-1\geq 2n-2i+1$. By Lemma~\ref{3l}, $e_1\geq n/2$. So, $\D+(d(w_1)+d(w_2)-1)+e_1\geq 3n-3i+2+n/2 \geq 3n$, a contradiction to~\eqref{e1}. Thus we can place all
 $v_1,\ldots,v_s$ on the corresponding $w_1,\ldots,w_s$.

 Next, we show that for every $i=1,\ldots,s$,
\begin{equation}\label{eq4}
 \mbox{the  number of edges incident to  vertices in $W_i=\{w_1,\ldots,w_i\}$ is at least $2i+1$.}
\end{equation}

 By Lemma~\ref{5l},~\eqref{eq4} holds for $i=1$.
Suppose~\eqref{eq4} holds for  some $i\leq s-1$. If  $w_{i+1}$ is not adjacent to $W_i$, then~\eqref{eq4} holds for   $i'=i+1$.
Otherwise, by the rules, $W_{i}\cup N(W_i)\supseteq V_2$ and the total number of edges incident to at least one vertex in $W_{i+1}$
is at least $n-(i+1)\geq n-s\geq 5n/6 \geq 2(i+1)+1$. This proves~\eqref{eq4}.

By~\eqref{eq4}, for $\mathbf{G}'=\mathbf{G}-H-W_s$, $|E(\mathbf{G}')|\leq |E(\mathbf{G})|-(s-1)-(2s+1)=|E(\mathbf{G})|-3s$. Then, $\mathbf{G}'$ does not pack, because $\mathbf{G}$ does not pack, and a packing of $\mathbf{G}'$ would extend to $\mathbf{G}$. By the minimality of $\mathbf{G}$, this yields ~\eqref{e0} does not hold. Then there exists some vertex $x$ such that $d_{j} (x)\geq n-s-1$ for some $j = 1,2,3$. Hence $\D \geq n - s - 5$. 

Now, we wish to say more about $H$. First, $H$ cannot be a single vertex by Lemma \ref{3l}. 
Suppose $H=K_2$. By Lemma~\ref{5l},  $d(w_1)\geq 3$. By~\eqref{dec243}, $d(w_2)\geq 2$.  In this case, the triple $\mathbf{G}'=\mathbf{G}-H-w_1-w_2$ 
has at most $e_1+e_2+e_3-6$ edges. So by~\eqref{e2} and the minimality of $\mathbf{G}$, triple $\mathbf{G}'$ packs,
 and this packing extends to $\mathbf{G}$ by placing $v_1$ on $w_1$ and $v_2$ on $w_2$.  Therefore, $s \geq 3$ and the average degree of $H$ is at least $\frac{4}{3}$. 
 In fact, since $H$ was the smallest tree component, all of $G_1$ has average degree at least $4/3.$ Thus,
\[ \D+\sum_{v\in V_1}d(v)\geq (n - s - 5) + \frac{4}{3}n = 2n + \frac{n}{3} - s - 5 \geq 2n + \frac{n}{3} - \frac{n}{6} - 1 > 2n,\]
contradicting our assumption.\QED

The next lemma uses Lemma~\ref{6l} and its proof is similar.

\begin{lem}\label{7l}
Every white tree-component in $G_1$ has at least $C/3$ vertices.
\end{lem}

\noindent\textbf{Proof:}  Suppose $T$ is a smallest white tree-component in $G_1$ and  $s:=|V(T)|\leq C/3$.
By Lemma~\ref{5l}, $G_{2}$ has a vertex $w$ of degree at least $3$.
If $T$ contains a vertex $v\notin N(w)$, then let $v_1=v$ and $w_1=w$. Otherwise, let $v_1$ be any
vertex of $T$ and $w_1$ be any non-neighbor of $v_1$ in $G_{2}$ (such $w_1$ exists by~\eqref{e2}).
Now we repeat some arguments from the proof of Lemma~\ref{6l}.

Consider $T$ as a rooted tree with root $v_1$, so that each $x\in V(T)-v_1$ has a unique parent in $T$.
Order the vertices of $T$: $v_1,\ldots,v_s$ in the Breadth-First-Order.
We will consecutively place all vertices of $T$ on vertices in $V_{2}$. We start by sending $v_1$ to $w_1$. For every $i=2,\ldots, s$,
if possible, we send $v_i$ to a vertex  $w_i\in V_{2}$ not adjacent to the image $w_{i'}$ of any $v_{i'}$ with $i'<i$.  If this is not possible, then just send $v_{i}$ to any
nonoccupied non-neighbor of the image $w_j$ of its parent $v_j$.

If we cannot choose a vertex to place some $v_{i}$, then each vertex of $V_{2}$ either
is a neighbor of both $v_i$ and $w_j$, where $v_j$ is the parent  of $v_i$, or is occupied by one of $v_1,\ldots,v_{i-1}$.
Thus $d_2(w_j)+d_3(v_i)+i-1\geq n$. Since $d_2(w_j)+d_3(v_i)+i-1\leq\D+4+C/3-1$, this
contradicts~\eqref{e2}.
  Thus we can place all
 $v_1,\ldots,v_s$ on some $w_1,\ldots,w_s$.

Let $W_i=\{w_1,\ldots,w_i\}$. If $d(w_1)\geq 3$, then~\eqref{eq4} holds for $i=1$. So we show that~\eqref{eq4} holds for each $i\leq s$
exactly as in the proof
of Lemma~\ref{6l}. In this case, for $\mathbf{G}'=\mathbf{G}-T-W_s$, $|E(\mathbf{G}')|\leq |E(\mathbf{G})|-(s-1)-(2s+1)=|E(\mathbf{G})|-3s$.
 If $d(w_1)= 2$, then $w$ (and each vertex of degree at least 3 in $V_2$) is adjacent to each vertex in $T$ and, in addition, we have
 an analog of~\eqref{eq4} with $2i$ in place of $2i+1$. So again, $|E(\mathbf{G}')|\leq |E(\mathbf{G})|-3s$.
By the choice of $\mathbf{G}$, the triple $\mathbf{G}'$ does not pack.
By the minimality of $\mathbf{G}$, this yields that~\eqref{e0} does not hold. Then $\D\geq n-s-5$, contradicting~\eqref{e2}.
\QED

\begin{claim}\label{12c}
For $i\in\{1,2\}$ and $u\in V_i$ there are at least $\frac{2C-16}{3}$ vertices in $V_i-N_i[u]$ of degree at most $3$.
\end{claim}

\noindent
\textbf{Proof:} We will use two cases.

\noindent \textbf{Case 1:}  $i=1$.  By~\eqref{dec243}, $\sum _{w\in V_2} d(w)\geq 2n$. So since
 $\D \geq d_1(u)$, we have \[\sum_{v\in{V_1-N_1[u]}}d(v)+4d_1(u)\leq \sum_{v\in{V_1-N_1[u]}}d(v)+\sum_{v\in{N_1[u]}}d(v)+2d_1(u) \leq 4n-2C.\]
 Therefore,
$\sum_{v\in{V_1-N_1[u]}}d(v)\leq  4(|V_1|-|N_1[u]|) + 4 - 2C.$

\noindent \textbf{Case 2:} $i=2$. Since $\D\geq d_2(u)$,
\begin{align*}
\sum_{v\in{V_2-N_2[u]}}d(v)+4d_2(u)&\leq \sum_{v\in{V_2-N_2[u]}}d(v)+3d(u)+d_2(u)\\
&\leq\sum_{v\in{V_2-N_2[u]}}d(v)+\sum_{v\in{N_2[u]}}d(v)+d_2(u)\\
&\leq 4n+12-2C,
\end{align*}
where $\D+\sum_{v\in V_2} d(v)\leq 4n+12-2C$ by Lemma \ref{6l}. Hence,
\[\sum_{v\in{V_2-N_2[u]}}d(v)\leq 4\left(\left|V_2\right|-\left|N_2[u]\right|\right)+16 - 2C.\]

Thus, in both cases, \[\sum_{v\in{V_i-N_i[u]}}d(v)\leq  4(|V_i|-|N_i[u]|) + 16 - 2C,\]
and the average degree of vertices in $V_i-N_i[u]$ is less than four. Since every vertex has positive degree, $V_i-N_i[u]$ contains at least $\frac{2C-16}{3}$ vertices of degree strictly less than 4.
\QED

For $i\in \{1,2\}$ and every $v\in V_i$, define the \emph{shared degree} of $v$, $\sd(v)$, as follows. 
If $d_i(v)<15$, then $\sd_i(v):= d_i(v)+\frac{2}{3}|\{x\in N_i(v)\,:\,d_i(x)\geq 15\}$ and
$\sd(v):=\sd_i(v)+d_3(v)$.
If $d_i(v)\geq 15$, then $\sd_i(v):= d_i(v)-\frac{2}{3}|\{x\in N_i(v)\,:\,d_i(x)< 15\}$ and  $\sd(v):=\sd_i(v)+d_3(v)$.
By definition, (a) $\sum_{v\in V_i}\sd_i(v)=2e_i$ and $\sum_{v\in V_i}\sd(v)=2e_i+e_3$, 
(b) $\sd(v)\geq d(v)$ if
$d_i(v)<15$, 
 (c) $\sd(v)\geq d(v)/3\geq 5$ if $d_i(v)\geq 15$, and (d) $3\sd(v)$ is an integer for every $v\in V_i$.

\begin{claim}\label{14c}
For $i\in\{1,2\}$ and $u\in V_i$, there is a vertex $v\in V_{3-i}-N[u]$ of shared degree at most $4$.
\end{claim}

\noindent
\textbf{Proof:} Let $S=V_{3-i}-N(u)$ and $s=|S|$. Suppose that $\sd(v)>4$ for every $v\in S$. Then
by the property (d) of shared degrees,
$\sum_{w\in S}\sd(w)\geq \frac{13}{3}s$. By Lemma~\ref{3l}
and properties (b) and (c) of shared degrees,
$\sum_{x\in V_{3-i}-S}\sd_{3-i}(x)\geq n-s$ and, since each vertex in $V_{3-i} - S$ is also a yellow neighbor of $u$, we have that $\sum_{x\in V_{3-i}-S}\sd (x)\geq 2(n-s)$.  Combining these two sums, we see that $2e_{3-i} + e_{3} = \sum_{x \in V_{3-i}} \sd (x) \geq \frac{13}{3} s+ 2(n-s)$.

If $i=1$, then by Lemma~\ref{7l}, $e_i=e_1\geq n(1-\frac{3}{C})$.
If $i=2$, then $\sum_{x\in V_i-u}d(x)\geq 2n-2$.  In both cases the yellow neighbors of $u$ were not included in the sum, so we have that
\[
\sum_{x \in V_i} d(x)\geq 2n \left( 1 - \frac{3}{C} \right) + (n-s).
\]
By definition, $\D\geq (d_3(u)-4)+\Delta_{3-i}\geq n-s-3$.
These inequalities and property (a) of shared degrees yield,
\begin{align*}
 2(e_1+e_2+e_3+\D) &\geq 2n \left( 1 - \frac{3}{C} \right) + (n-s) + 2(n-s)+\frac{13}{3}s + 2(n-s-3) \\
&=\left(7-\frac{6}{C}\right) n - \frac{2}{3} s - 6 > 6n - 6.
\end{align*}
By~\eqref{e1}, this is at most $6n-2C$, a contradiction.\QED

\begin{lem}\label{8l}
Let $F:= \displaystyle\sqrt{\frac{C}{11}}= 195$. Then $D_{1}, D_{2} \geq F$.
\end{lem}

\noindent
\textbf{Proof:}  Suppose that $D_1\leq D_2$ and  $D_1 < F  = \sqrt{C/11}$; the proof for $D_2$ is similar. 
 By Theorem~\ref{List S-S product}, $D_{2}F + D_{2} \geq D_2D_1+\max\{D_1,D_2\} \geq n/2$, so $D_{2} \geq n/ (2F + 2) $.  Consider a vertex $w \in V_{2}$ of maximum degree. By the choice, $d(w) \geq D_{2}$. By \eqref{e2}, $d_2(w) < n - C + 2$.  By Claim~\ref{14c}, $V_1$ contains a non-neighbor $v$ of $w$ with $\sd(v)\leq 4$.
 In particular, by the definition of shared degree, $d(v)\leq 4$.
 Let  $N_1(v) := \{v_1,\ldots,v_s\}$. 
 We wish to find an independent set  $\{w_1,\ldots,w_s\}\subset V_{2} - N_{2}[w]$ such that each $w_i$ has degree at most $3$ and is not adjacent to $v_i$.

By Claim~\ref{12c}, at least $\frac{2C-16}{3}$ vertices in $V_{2} - N_{2}[w]$ have degree at most $3$.  At most $F-1$ of them are adjacent to $v_1$. 
So, we can choose $w_1\in V_{2} - N_{2}[w]-N(v_1)$ with $d(w_1)\leq 3$. Continuing in this way for $j=2,\ldots,s$, at least $\frac{2C-16}{3} - 4 (j-1)$ vertices in $V_{2} - N_{2}[w]-\bigcup_{i=1}^{j-1}N[w_i]$ have degree at most $3$.  Again, at most $F-1$ of them are adjacent to $v_j$. Since $s\leq 4$ and
 $\frac{2C}{3} - 5 -4(s-1)-F \geq \frac{2C-16}{3} -17- F > 0$, we can choose $w_j\in V_{2} - N_{2}[w]-\bigcup_{i=1}^{j-1}N[w_i]-N(v_j)$ with $d(w_j)\leq 3$. 

We now create a new graph triple $\mathbf{G}' = (G_1', G_2', G_3')$ by removing $\{w,v,w_1,\ldots,w_s,v_1,\ldots,v_s\}$ and adding new yellow edges between $N_1(v_i)$ and $N_2(w_i)$ for each $1 \leq i \leq s$ and then deleting the matched pairs.  Through this process, since the set $\{w_1,\ldots,w_s,w\}$ is independent, we have removed at least 
$d(v) + d(w) + \sum_{i = 1}^{s} (d_1 (v_i) -1  + d_2(w_i))-|E(G_1[N_1(v)])|$ edges, 
and added at most $3 \sum_{i = 1}^{s} (d_1 (v_i) -1) - 2|E(G_{1}[N_{1}(v)])|$ edges. We have increased 
$\D$ by at most $\max\{\max_i (d_1(v_i) - 1 ), \max_j d_2(w_j)\}\leq F-1$.  Thus, we have
\[ \partial (\mathbf{G},\mathbf{G}')\geq d(v) + d(w) + \sum_{i = 1}^{s} d_2(w_i) - 2\sum_{i = 1}^{s} (d_1 (v_i) -1) - F + |E(G_{1}[N_{1}(v)])| + 1, \]
and therefore
\begin{equation}\label{dec28}
\partial (\mathbf{G},\mathbf{G}')\geq  d(w)- 2\sum_{i = 1}^{s} (d_1 (v_i) -1)-F.
\end{equation}
If $s\leq 2$, then $\sum_{i = 1}^{s} (d_1 (v_i) -1)\leq 2F-2$. If $s=3$, then since $\sd(v)\leq 4$, at least two neighbors of $v$ have degree less than $15$,
so in this case $\sum_{i = 1}^{s} (d_1 (v_i) -1)\leq 2\cdot 13+F-1=25+F\leq 2F-2$.
If $s=4$, then since $\sd(v)\leq 4$, all $4$ neighbors of $v$ have degree less than $15$. So in this case $\sum_{i = 1}^{s} (d_1 (v_i) -1)\leq 4\cdot 13\leq 2F-2$.
So since $d(w) \geq D_{2} \geq \frac{n}{2(F+1)}\geq \frac{C}{2F+2}$, by~\eqref{dec28} and the definitions of $C$ and $F$, 
\[\partial (\mathbf{G},\mathbf{G}')\geq \frac{C}{2F+2}-2(2F-2)-F=\frac{C}{2F+2}-5F+4\geq 15\geq 3(s+1).\]
It follows that~\eqref{e1} holds for $\mathbf{G}'$.
Also by above, $\D'-\D\leq  F-1$. Thus by~\eqref{e2},
  \[\D'\leq \D+F-1\leq n+2-C+F-1=(n'+s+1)+1-C+F<n'-5,\]
and~\eqref{e0} holds for $\mathbf{G}'$.
So   $\mathbf{G}'$ packs by the minimality of $\mathbf{G}$, and then $\mathbf{G}$ also packs, a contradiction.  
  \QED

\begin{lem}\label{10l}
Let $K:=\displaystyle\frac{F}{13}=15$. Let $i\in \{1,2\}$ and  $v \in V_{i} $ with $d(v)=t \leq 4$ be
 not adjacent to some vertex $w\in V_{3-i}$ of degree at least~$F$.\\
(a)  Then  $v$ has a neighbor in $V_i$ of degree at least $\frac{13K}{3t+1}$.\\
(b) Moreover, if $2\leq t\leq 3$ and $v$ has $t-1$ neighbors of degree at most $2$, then
$v$ has a neighbor in $V_i$ of degree at least $\frac{13K}{5}$.
\end{lem}

\noindent
\textbf{Proof:} Suppose Statement (a) of the lemma fails for $i=1$ (the proof for $i=2$ is the same). This means that
for  a vertex $v\in V_1$ of degree  $t$ in $\mathbf{G}$,  all of its neighbors  in $V_1$  have degree less than $\frac{13K}{3t+1}$ and
some  non-neighbor $w\in V_2$ of $v$ has $d(w) \geq F$.
 Let $N_1(v):=
 \{v_1,\ldots,v_s\}$.   By definition, $s\leq t\leq 4$.
  We wish to find an independent set  $\{w_1,\ldots,w_s\}\subset V_{2} - N_{2}[w]$ 
 such that each  $w_i$ has degree at most $3$ and is not adjacent to $v_i$.

By Claim~\ref{12c}, at least $\frac{2C-16}{3}$ vertices in $V_{2} - N_{2}[w]$ have degree at most $3$.  
Less than $\frac{13K}{3t+1}-1$  of them are adjacent to $v_1$. So, we can choose $w_1\in V_{2} - N_{2}[w]-N(v_1)$ with
$d(w_1)\leq 3$. Continuing in this way for $j=2,\ldots,s$,
 at least $\frac{2C-16}{3}-4(j-1)$ vertices in $V_{2} - N_{2}[w]-\bigcup_{i=1}^{j-1}N[w_i]$ have degree at most $3$. 
Again less than $\frac{13K}{3t+1}-1$ of them are adjacent to $v_j$. Since $\frac{2C-16}{3}-4s-\frac{13K}{3t+1}\geq\frac{2C-16}{3}-16-\frac{13K}{3t+1}>0$, 
we can choose $w_j\in V_{2} - N_{2}[w]-\bigcup_{i=1}^{j-1}N[w_i]-N(v_j)$ with
$d(w_j)\leq 3$. 

Finally, we can map $v$ to $w$, vertices $v_1,\ldots,v_s$ to $w_1,\ldots,w_s$, respectively, delete the matched pairs, and for each  pair $\{v_i, w_i\}$, introduce yellow edges between the remaining vertices of $N_1(v_i)$ and $N_2 (w_i)$. 
This creates a new graph triple $\mathbf{G}' = (G_1', G_2', G_3')$.  
During this process, we have deleted at least $d(w) + d(v) $ edges,  added in strictly less than $3 s (\frac{13K}{3t+1}-1)$ new yellow edges, and increased $\D$ by at most $\max\{3,\max_i\{d_1(v_i)-1\}\}\leq \frac{13K}{3t+1}-1$.  
Therefore since $F=13K$,
\begin{align}\label{dec29}
\partial (\mathbf{G},\mathbf{G}') &> d(v) + d(w)  - (3s+1) \left( \frac{13K}{3t+1} - 1 \right) \nonumber \\
&\geq s + d(w) -  13K+(3s+1) \\
&\geq  F - 13 K+(4s+1) \geq 3s+2. \nonumber
\end{align}
Now, we need $\partial (\mathbf{G},\mathbf{G}')\geq 3s+3$ but since we added \emph{strictly} less than $3 s (\frac{13K}{3t+1}-1)$ yellow edges, we have a strict inequality which, in combination with the fact that both $\partial (\mathbf{G},\mathbf{G}')$ and $3s+2$ are integers, in fact gives $\partial (\mathbf{G},\mathbf{G}')\geq 3s+3$. 
Since $\partial (\mathbf{G},\mathbf{G}')$ is sufficiently large and $\mathbf{G}$ is a minimal counterexample, $\mathbf{G}'$ packs unless \eqref{e0} is violated.  However, by \eqref{e2}, this violation would have to occur at some vertex in some $N_1(v_i)$ or $N_2(w_i)$ but the degrees of these vertices only increase by at most $3$ or 
$ (\frac{13K}{3t+1}-1)<4K$, neither of which could get us to have a vertex of degree $(n-s-1)-2\geq n-7 $. Hence, $\mathbf{G}'$ packs and this packing extends to a packing of $\mathbf{G}$, as we constructed above.  This proves (a).

To prove (b), we repeat the argument of (a) with $\frac{13K}{5}$ in place of $\frac{13K}{3t+1}$ until we count the number of added
yellow edges. We have added less than $3\left((s-1)+\frac{13K}{5}\right)$ edges and increased $\D$ by 
 at most $ \frac{13K}{5}-1$.  So, instead of~\eqref{dec29}, we will have
\begin{align*}
\partial (\mathbf{G},\mathbf{G}') &> d(v) + d(w)  - 3(s-1)-4 \left( \frac{13K}{5} - 1 \right) \\
&\geq s+13K-3(s-1)-\frac{4\cdot 13K}{5}+4 \\
&=\frac{13K}{5}-2s+7>3s+3.
\end{align*}
Then again we simply repeat the last paragraph of the proof of (a).
\QED

%======================================================================
%Section:  At Most One Vertex in $V_1$ has two 1-neighbors
%======================================================================
% 
%
\section{At Most One Vertex in $V_1$ is a donor}\label{sec:1-neighbors}
%
%
%======================================================================

Recall that  by Lemma~\ref{4l} we assumed (see~\eqref{dec243}) that $V_{2}$ has no vertices of degree $1$. 
A {\em donor} is a vertex in $V_1$ adjacent to at least two vertices of degree $1$.
The goal of this section is to prove that $V_1$ contains at most one donor.

\begin{lem}\label{11l}
Suppose  $V_{1}$ contains   donors $v$ and $v'$.  If $w \in V_{2}$ with $d(w) = 2$, then
$N(w)\subset V_2$ and $d(w') \geq 2K$ for each $w' \in N(w)$.
\end{lem}

\noindent
\textbf{Proof: }  Suppose the lemma fails for some  $w \in V_{2}$ with $d(w) = 2$.
Let $x,y \in V_1$ be degree one neighbors of $v$ and let $x',y' \in V_1$ be degree one neighbors of $v'$.  By Lemma~\ref{10l}, $d(v), d(v') \geq 3K$.

\textbf{Case 1: } $N(w) = \{ w_1, w_2 \} \subset V_2.$ By symmetry, assume  $d(w_2) < 2K $.  
Begin by mapping $x$ and $y$ to $w_1$ and $w_2$ , respectively, and adding new yellow edges from $N_2(w_1) \cup N_2(w_2) - \{w\}$ to $v$.  Since $v$ is the only neighbor of $x$ and $y$, this assignment is permitted and adding the yellow edges ensures that any permissible extension of the mapping will not violate the packing property.  
After mapping $x$ and $y$,  $w$ is adjacent only to $v$ and so $v'$ may be mapped to w.  This in turn causes $x'$ and $y'$ to be newly isolated vertices.  After removing these 3 pairs of vertices and adding the yellow edges, let $z \in V_2 - \{w, w_1, w_2 \}$ be the vertex of $V_2$ of highest degree and map $x'$ to $z$.

We now have a new graph triple $\mathbf{G}' := (G_1', G_2', G_3 ')$.  Note that $\Delta_1',\Delta_2'\leq n'-2$ since \eqref{e2} holds for $\mathbf{G}$ so that \eqref{e0} is only violated if $d_3'(v)=n-4$. However, 
\[d_3'(v)\leq (d_3(v)+d_2(w_1))+d_2(w_2)\leq (\D+4)+2K \leq n-C+6+2K<n-4,\] 
so \eqref{e0} is satisfied for $\mathbf{G}'$ as well. Now, we will consider $\partial (\mathbf{G},\mathbf{G}')$.  In particular, we have deleted at least $d(w_1) + d(w_2) - \| w_1, w_2\|$ edges adjacent to $w_1$ and $w_2$ and exactly $2$ edges adjacent to $x$ and $y$.  We then added at most $(d_2 (w_1) -1) + (d_2 (w_2) - 1) -  |N_2 (w_1) \cap N_2 (w_2) - \{w\}| -  2\| w_1, w_2\|$ yellow edges.  Finally, we deleted at least $d(v') - 1 - \|v', \{w_1, w_2 \} \|$ edges adjacent to $v$ and at least $d(z) - \max \{ 0, \| z, \{w_1, w_2 \}\| - 1\}$ edges adjacent to $z$.  To see this, note that if $\| z, \{w_1, w_2 \}\| \neq 0$, then we save one additional edge, since $vz$ must now be a yellow edge in the modified graph (either $vz \in E_{3}$ and we didn't need to add it to begin with, or it was added and the degree of $z$ grew by one before we deleted it).  In any event, $|N_2 (w_1) \cap N_2 (w_2) - \{w\}| - \max \{ 0, \| z, \{w_1, w_2 \}\| - 1\} \geq 0$.  Thus,
\[d(w_1) + d(w_2) + \|w_1, w_2\| \geq d_2(w_1) + d_2(w_2) + \|v', \{w_1, w_2 \} \| . \]
 Therefore, the total change in the number of edges is:
\begin{equation}\label{20l}
e(\mathbf{G}) - e(\mathbf{G}') \geq d(v')  + d(z) +1.
\end{equation}

Next, consider the difference $\D - \D'$.  If $\D - \D' \geq -1$, then $\partial (\mathbf{G}, \mathbf{G}') \geq d(v') + d(z) \geq 12$ and $\mathbf{G}'$ packs by the inductive assumption.  If $\D - \D' \leq -2$, then we must have that $\D' = d'_{3} (v) + \Delta_2' - 4$.  In particular, since $d(z) \geq \Delta_{2}'$, $\Delta_{2} \geq d_{2}(w_{1})$, and $d_{3}(v) - d_{3}'(v) \geq 2 - d_{2}(w_{1}) - d_{2} (w_{2})$, 
\[ \D - \D' \geq 2 - d_{2}(w_{1}) - d_{2} (w_{2}) + d_{2}(w_{1}) - d(z) = 2 - d_{2}(w_{2}) - d(z).\]
Combining this with \eqref{20l} , we see that
\[ \partial (\mathbf{G},\mathbf{G}') \geq (d(v')  + d(z) +1) + (2 - d_{2}(w_{2}) - d(z)) = d(v') - d_{2}(w_{2}) + 3.\]
Since $d(v') \geq 3K$ and $d(w_2) \leq 2K $, we have $\partial (\mathbf{G},\mathbf{G}') \geq 12$.  By the minimality of $\mathbf{G}$, we conclude that $\mathbf{G}'$ packs. 
And  we can extend any packing of $\mathbf{G}'$ to a packing of $\mathbf{G}$.

\textbf{Case 2: } $N_2(w) = \{w'\}$.  This case follows in a similar fashion to Case 1.  Since $d_3(w) = 1$, we may assume that $v' \notin N(w)$.  We begin by mapping $x$ to $w'$ and adding new yellow edges from $v$ to $N_2(w') - w$.  We then map $v'$ to $w$ and choose a remaining vertex $z \in V_2$ of maximum degree to have $x'$ map to z.  
Then we delete the matched pairs.
This process creates a new graph triple $\mathbf{G}'' := (G_1'', G_2'', G_3'')$. Again, the only way \eqref{e0} is violated is if $d_3'(v)=n-3$, but this is not the case, since \[d_3('(v)\leq d_3(v)+d_2(w')\leq \D+4\leq n+6-C<n-3.\]

During this process, we removed $d(w')$ edges adjacent to $w'$, one edge adjacent to $x$, one yellow edge adjacent to $w$, at most $d(v') - 1 - \|v', w'\|$ edges adjacent to $v'$, and $d(z) - \| w', z \|$ edges adjacent to $z$.  We have added in $d_2(w') - 1 - \|w', z \|$ new yellow edges.  Since $d(w') \geq d_2(w') + \| v', w\|$, we see that:
\[e(\mathbf{G}) - e(\mathbf{G}'') \geq  d(v')  + d(z) + 2. \]

As in Case 1, if $\D - \D' \geq -1$, then $\partial (\mathbf{G}, \mathbf{G}') \geq d(v') + d(z) \geq 12$ and $\mathbf{G}''$ packs by the inductive assumption.  If $\D - \D' \leq -2$, then we must have that $\D' = d'_{3} (v) + \Delta_2' - 4$.  Since $d(z) \geq \Delta_{2}'$, $\Delta_{2} \geq d_{2}(w')$, and $d_{3}(v) - d_{3}'(v) \geq 1 - d_{2}(w')$, we must have that
\[ \D - \D' \geq 1 - d_{2}(w') + d_{2}(w') - d(z) = 1 - d(z).\]
Thus, 
\[ \partial (\mathbf{G},\mathbf{G}') \geq (d(v')  + d(z) +1) + (1 - d(z)) \geq d(v') + 2 \geq 9.\]
By the minimality of $\mathbf{G}$, triple $\mathbf{G}'$ has a packing, which we can extend to a packing of $\mathbf{G}$. 

\begin{flushright}
\QED
\end{flushright}

\begin{cor}\label{cor1}
Suppose  $V_{1}$ contains   donors $v$ and $v'$.
 Then $\displaystyle 2e_{2} + e_{3} = \sum_{v \in V_{2}} d(v) \geq 3n$.
\end{cor}

\noindent
\textbf{Proof:}  Consider the following discharging.  For each vertex $v \in V_{2}$, assign $v$ charge $d(v)$.  The total charge allocated is $\sum_{v \in V_{2}} d(v) = 2e_{2} + e_{3}$.  Now, each vertex of degree at least $6$ will give charge $\frac{1}{2}$ to each neighbor and save $d(v) / 2 \geq 3$ for itself.  By Lemma~\ref{11l}, each vertex of degree $2$ is adjacent to two vertices in $V_{2}$ with degree at least $2K  \geq 30$.  Thus, after discharging each vertex has charge at least $3$.  So the total charge is at least $3n$ and  $2e_{2} + e_{3} \geq 3n$, as needed. \QED

\begin{rem}\label{r1}
Suppose  $V_{1}$ contains   donors $v$ and $v'$.
  If $w \in V_{2}$ with $d(w) = 3$
and $v'w\notin E(\mathbf{G})$, then $w$ has a neighbor in $V_2$ of degree at least
$K+1$.
\end{rem}

\noindent
\textbf{Proof: } If $w$ has no yellow neighbors, this follows from Lemma~\ref{10l}. Otherwise, suppose the remark fails for some $w \in V_{2}$ with $d(w) = 3$.
 Then  each of the neighbor(s)  $w_1$ and $w_2$ (if it exists) 
of $w$ in $V_2$  has degree at most $K$.  Map $w$ to $v'$ and map two degree one neighbors of $v$ to $w_1$ and $w_2$.  Next, form a new graph triple $\mathbf{G}'$ by adding new yellow edges from $v$ to $W := N_2(w_1) \cup N_2(w_2) - \{ w,w_1,w_2\}$ and deleting the previously matched pairs.  We have deleted at least $d(v') +  2 + d_2(w_1) + d_2(w_2)-\|w_1,w_2\|$ edges
  and added $|W|$ new yellow edges.  We have increased $\D$ by at most $|W|$.  Since $d(w_1) + d(w_2)-\|w_1,w_2\|-1 \geq |W|$ 
  (in fact, it is at least $|W| + 1$ if $w_2$ exists),   $\partial (\mathbf{G},\mathbf{G}')\geq d(v') + 3 - |W|$.  Now $|W| \leq 2K-2$ and $d(v') \geq 3K$, so that $\partial (\mathbf{G},\mathbf{G}') \geq 12$.  In particular, by the minimality of $\mathbf{G}$,
$\mathbf{G}'$ has a packing, and it extends to a packing of $\mathbf{G}$, a contradiction. \QED

\begin{lem}\label{12l}
Suppose  $V_{1}$ contains   donors $v$ and $v'$.
Then $\D\leq  \frac{9n}{4K} $.  
\end{lem}

\noindent
\textbf{Proof: } Suppose $\D >  \frac{9n}{4K} $. By Lemma~\ref{7l}, $e_1\geq n(1-3/C)$. 

Consider the following discharging on $V_2\cup E_3$. The initial charge, $ch(v)$, of every $v\in V_2$ is $d(v)$ and of every edge in $E_3$ is $1$. 
The total sum of charges, $ch(w)$, over  $w\in V_2\cup E_3$ is  $2(e_2+e_3)$. We use two rules.

(R1) Each vertex $w\in V_2$ of degree at least $5$ gives to every neighbor in $V_2$ charge $\frac{d(w)-4}{d(w)}$. 

(R2) Each edge in $E_3$ gives charge $1$ to its end in $V_2$.

Let $ch^*(w)$ denote the new charge of $w\in V_2\cup E_3$. By (R2), $ch^*(w)=0$ for every $w\in E_3$.
 By (R1), if $w\in V_2$ and $d(w)\geq 4$,
then $ch^*(w)\geq 4$. If $d(w)=3$ then by (R1), (R2) and Lemma~\ref{10l}, $ch^*(w)\geq 3+ (1-\frac{4}{K})$. 
If $d(w)=2$ then by Lemmas~\ref{10l} and~\ref{11l}, 
	\[ch^*(w)\geq 2+ 2(1-\frac{2}{K})=4-\frac{4}{K}.\]
Since the total sum of charges did not change, we conclude that
\[2( e_2+e_3)=\sum_{w\in V_2}ch^*(w)\geq 4n\left(1-\frac{1}{K}\right).\]
It follows that
\begin{align*}
e_1+e_2+e_3+\D &\geq n\left(1-\frac{3}{C}\right)+n\left( 2 - \frac{2}{K} \right)+ n \left(  \frac{9}{4K} \right) \\
&\geq 3n+n\left(  - \frac{3}{C} +\frac{1}{4K}\right).
\end{align*}
Since ${4K}\leq \frac{C}{3}$, this
 contradicts~\eqref{e1}. \QED

For $v\in V_1$, let $L(v)$ be the set of neighbors of $v$ of degree $1$.

\begin{lem}\label{13l}
Suppose  $V_{1}$ contains   donors $v$ and $v'$.
Then $|L(x)| \leq d(x) / 2$ for every $x\in V_1$. 
\end{lem}

\noindent
\textbf{Proof:} Suppose $x\in V_1$, $\ell = |L(x)| > d(x) / 2$ and $L(x) = \{ x_{1}, \ldots, x_{\ell} \}$.  By Lemma~\ref{10l}, $d(x) \geq K$.
Thus, $x$ is a donor, so we may assume $x=v$.

\textbf{Case 1:} \emph{There is a vertex $w \in V_{2} - N_{3}(v)$ with $d_{2}(w) \leq 2$.}
Let $w_{1}$ be a white neighbor of $w$ and, if it exists, let $w_{2}$ be the other white neighbor of $w$.  We wish to find a vertex in $V_{2}-\{w,w_1,w_2\}$ 
with low degree that is  adjacent to none of $w_{1}$, $w_{2}$, or $v'$.  By Lemma~\ref{12l} and since $K = 15$, we have $\D \leq \frac{9n}{4K}=\frac{3n}{20}$.  
By definition, $d_2(w_1)+(d_3(v')-4)\leq \D$.
Therefore, 
\[|V_{2} - N[ \{w_{1}, w_{2}, v' \}]|\geq
(n-3) -  \D - (\D + 4) \geq \frac{14n}{20} - 7 \geq \frac{n}{2}.\]
  Since $\sum_{w \in V_{2}} d(v) < 4n$ by Lemma~\ref{6l} and \eqref{e1}, the average degree of the vertices in $V_{2} - N[ \{w_{1}, w_{2}, v' \}]$ is less than $8$.
 So, there exists   a vertex $w'\in V_{2} - N[ \{w_{1}, w_{2}, v' \}]$ with $d(w')\leq 7$.

Construct a packing in the following way.  Since $\ell \geq \frac{13}{8} K  > 7$, we may send $x_{1}, \ldots, x_{d_{2}(w')}$ to the white neighbors of $w'$.  Send two degree $1$ neighbors of $v'$ to $w_{1}$ and $w_{2}$.  Finally, send $v$ to $w$ and $v'$ to $w'$. Let $\mathbf{G}'$ be obtained by deleting the matched pairs.
 Then $n - n' \leq 11$.  By Lemma~\ref{10l}, we have deleted at least $d(v) + d(v') - \| v, v' \| \geq \frac{13}{2}K - 1 \geq 36$ edges and~\eqref{e0} still holds, so $\mathbf{G}'$ packs.  This packing extends to a packing of $\mathbf{G}$, a contradiction.

\textbf{Case 2:} \emph{Every vertex $w \in V_{2} - N_{3}(v)$ has $d_{2}(w) \geq 3$.}
If there is a vertex $w \in V_{2}$ with $d(w) = 2$, then $N(w) \subset V_{2}$ by Lemma~\ref{11l} and we have Case~1.  So, $d(w) \geq 3$ for all $w \in V_{2}$.
 If every vertex in $X := V_{1} - N_{1}[v]-N_1[v']$ has degree at least $3$,  then
\begin{align}
\sum_{x \in V_{1}} d(x) + 2\D &= \sum_{x \in N_{1} (v)\cup N_1(v')} d(x) + \sum_{y \in X} d(y) + d(v)+d(v') + 2\D \nonumber \\
&\geq d_{1}(v) +d_1(v')+ 3(n - 2 - d_{1}(v)-d_1(v') ) +  d(v)+d(v') + 2\D  \label{13l:case2}  \\
&\geq 3n - 6. \nonumber
\end{align}
Since every vertex in $V_{2}$ has degree at least $3$, we get
\[ \sum_{x \in V} d(x) + 2\D \geq (3n -6) + 3n \geq 6n-6,\]
a contradiction to \eqref{e1}.  So there is a vertex $v_{0} \in V_{1} - N_{1}[v]-N_1[v']$ with $d(v_0) \leq 2$.

By Lemma~\ref{6l} and \eqref{e1}, $\sum_{v \in V_{2}} d(v) + \D \leq 4n - 2C + 12$ and so there are at least $2C + \D - 12$ vertices of degree $3$ in $V_{2}$.  Moreover, since $d_{3}(v) \leq \D + 4$, there is a vertex $w \in V_{2} - N_{3}(v)$ with $d(w) = 3$.  By Case~1, all neighbors of $w$ are white so let $\{w_{1}, w_{2}, w_{3} \}=N_2(w)$ with 
\begin{equation}\label{dec263}
d_{2} (w_{1}) \geq d_{2} (w_{2}) \geq d_{2} (w_3)\geq 3.
\end{equation}
  Similarly to Case~1, we wish to find a vertex in $V_{2}$ with low degree that is  adjacent to none of
 $w_{1}, w_{2}, w_{3}, v'$.  As in Case 1, we use $d_2(w_1)+(d_3(v')-4)\leq \D$. This yields that
  \[|V_{2} - N[ \{w_{1}, w_{2}, w_{3}, v' \}]|\geq (n-4) - 2 \D - (\D + 4) \geq \frac{11n}{20} - 8 \geq \frac{n}{2}.\]
    Since $\sum_{w \in V_{2}} d(v) < 4n$ by Lemma~\ref{6l} and \eqref{e1}, the average degree of $V_{2} - N[ \{w_{1}, w_{2}, w_{3}, v' \}]$ is less than $8$ and there exists a vertex $w'$ in this set with degree at most $7$.  
    \begin{figure}[h!]
\begin{center}
\includegraphics[width=.5\textwidth]{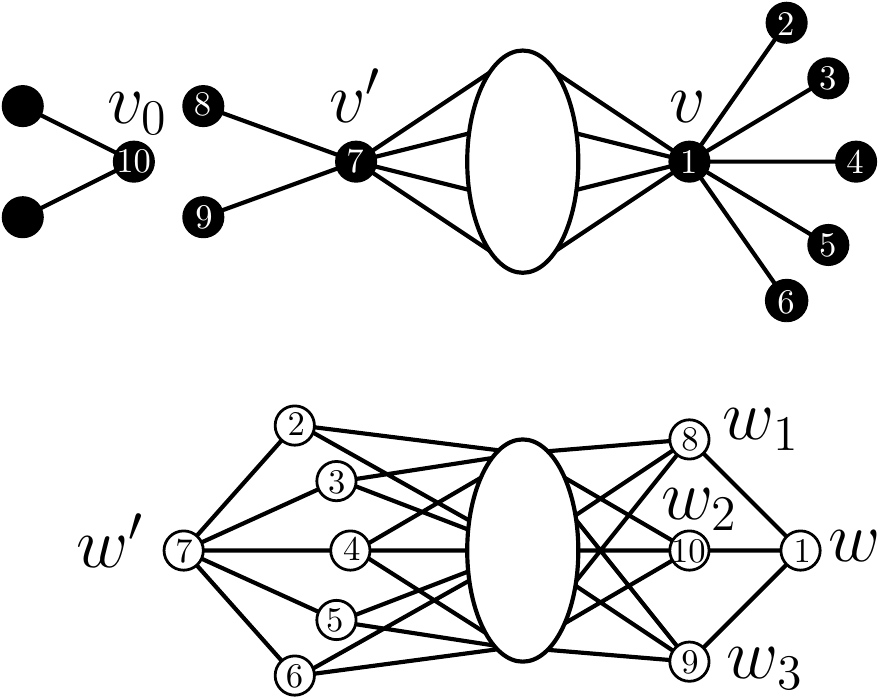}
\end{center}
\caption{Sketch of the packing used in Lemma~\ref{13l}}
\label{fig:13l}%
\end{figure}

Let $ j$ be the largest index such that $v_0w_j\notin E_3$ and $j\leq 3$.
Since $d(v_0)\leq 2$ and $v_0$ has a neighbor in $V_1$, $\| v_0,\{w_1,w_2,w_3\} \| \leq 1$. So, $j\geq 2$.

 Since $\ell \geq \frac{13}{8} K  > 7$, we may send $x_{1}, \ldots, x_{d_{2}(w')}$ to the white neighbors of $w'$.  Send two degree $1$ neighbors of $v'$ to 
 the vertices in $\{w_1,w_2,w_3\}-w_j$ and $v_{0}$ to $w_{3}$. 
   Send $v$ to $w$ and $v'$ to $w'$.  Finally,  add yellow edges between the white neighbors of $v_{0}$ and the white neighbors of $w_{j}$.  
   Delete the matched pairs.
   The resulting triple $\mathbf{G}'$ has order $n - 5  -  d_{2}(w')$.  We added at most $d_1(v_0)(d_{2}(w_{j})-1)\leq 2(d_{2}(w_{j})-1)$ yellow edges, and 
\begin{equation}\label{dec262}
 \D' \leq \D + \max\{2,d_{2} (w_{j})-1\}\leq 2\D-1.
\end{equation}   
  By Lemma~\ref{12l} and~\eqref{dec262}, \eqref{e0} holds. 
    The number of deleted edges is at least
\[ d_{2}(w') + d_{2} (w_{1}) + d_{2} (w_{2}) + d_{2} (w_{3}) -|E(G_2[\{w_1,w_2,w_3\})|+ d(v) + d(v') - \| v , v' \|+ d(v_0). \]
\begin{equation}\label{dec264}
\geq d_{2}(w') + d_{2} (w_{1}) + d_{2} (w_{2}) + d_{2} (w_{3}) -4+ d(v) + d(v')+d(v_0).
\end{equation}
{\bf Case 2.1:} $j=3$. Then by~\eqref{dec262}, the number of added yellow edges plus $\D'-\D$ is at most $3(d_2(w_3)-1)+\max\{3-d_2(w_3),0\}$.
Since $d_2(w_3)\geq 1$, by~\eqref{dec263}, this is at most $d_{2} (w_{1}) +d_2(w_2)+d_2(w_3)-1$. So by~\eqref{dec264} and because $d(w')\leq 7$,
\begin{equation}\label{dec265}
\partial (\mathbf{G},\mathbf{G}') \geq d_{2}(w') + d(v) + d(v')-2 \geq d_{2}(w') + \frac{13}{2} K - 2 \geq 3 (d_{2}(w') + 5). \end{equation}  
Therefore,  $\mathbf{G}'$ packs by the minimality of $\mathbf{G}$, and this packing extends to a packing of $\mathbf{G}$, a contradiction.

{\bf Case 2.2:} $j=2$. By the choice of $j$, this means $v_0w_3\in E_3$. Since $d(v_0)\leq 2$ and $v_0$ has a white neighbor,
$d(v_0)= 2$ and $d_1(v_0)=1$. It follows that we have added at most $d_{2}(w_{2})-1$ yellow edges, and so by~\eqref{dec264}, 
similarly to~\eqref{dec265},
we get
\[\partial (\mathbf{G},\mathbf{G}') \geq d_{2}(w')  + d_{2} (w_{3})+ d(v) + d(v')-2 \geq d_{2}(w') + \frac{13}{2} K - 2 \geq 3 (d_{2}(w') + 5),\]
which similarly yields a contradiction.
\QED

\begin{lem}\label{14l}
 $ V_{1}$  contains at most one donor.
\end{lem}

\noindent
\textbf{Proof: } 
Suppose  $v$ and $v'$ are donors in  $V_{1}$.  
Consider the following discharging. 

At start, we let $ch(v)=d(v)+\D + 4$,  $ch(v')=d(v')+\D + 4$, and $ch(u)=d(u)$ for each $u\in V(\mathbf{G})-v-v'$. 
By definition, the total sum of charges is $\sum_{v\in V(\mathbf{G})}d(v)+2\D+8=2F(\mathbf{G})+8$.
We redistribute charges 
according to the following rules.

(R1)  Each vertex $u$ not adjacent to $1$-vertices with $d(u) \geq 4$ gives to each neighbor
charge $\frac{d(u)-4}{d(u)}$ (and keeps $4$ for itself). 

(R2) Each vertex $x$ adjacent to $1$-vertices (it must be in $V_1$ and have degree at least $3K$) gives to each $z\in L(x)$ charge $\frac{4}{3}$ and to each $z'\in N(x)-L(x)$ charge $\frac{|N(x)-L(x)| - \frac{1}{3} |L(x)|-3}{|N(x)-L(x)|}$.

(R3) Each of $v,v'$, in addition, gives $1$ to each yellow neighbor.

\medskip
We will show that the resulting charge, $ch^*$, satisfies
\begin{equation}\label{5273}
ch^*(x)\geq \frac{7}{3} \quad\mbox{for each $x\in V_1$}\qquad\mbox{and} \qquad ch^*(y)\geq \frac{11}{3} \quad\mbox{for each $y \in V_2$}.
\end{equation}
This would mean that $\sum_{v\in V(\mathbf{G})}d(v)+2\D + 8 \geq \frac{7}{3} n + \frac{11}{3} n = 6n$, a contradiction to~\eqref{e1}.

If $d(u)=1$, then $u\in V_1$ and by (R2), $ch^*(u) = d(u) + \frac{4}{3} = \frac{7}{3}$, as claimed.
If $d(u) = 2$ and $u \in V_{1}$, then by  Lemma~\ref{10l}, $u$ has a neighbor $x$ with $d(x) \geq \left\lceil\frac{13K}{7}\right\rceil=28$. If $x$ has no neighbors of degree $1$, then by (R1) it gives to $u$ charge $\frac{d(x)-4}{d(x)}\geq 1-\frac{4}{28}>\frac{1}{3}$.
Otherwise, by (R2), it gives to $u$ charge $\frac{ |N(x)-L(x)| - \frac{1}{3} |L(x)|-3}{|N(x)-L(x)|}$. By Lemmas~\ref{13l} and~\ref{10l}, this is at least 
$1 - \frac{1}{3} - \frac{3}{|N(x)-L(x)|} \geq \frac{2}{3}-\frac{3}{28/2}> \frac{1}{3}$.
If $d(u)=2$ and $u\in V_2$, then by  Lemma~\ref{11l}, both neighbors of $u$ are in $V_2$, and each of them has degree at least $2K$. 
So by (R1), $ch^*(u) \geq 2+  2\frac{2K-4}{2K} = 4-\frac{4}{K} = 4-\frac{4}{15}>\frac{11}{3}$.

If $d(u)\geq 3$, $u\in V_1$ and $u$ has no neighbors of degree $1$, then either $u$ keeps all its original charge (when $d(u)\leq 4$) or keeps for itself charge $4$ by (R1). In both cases,~$ch^*(u)\geq 3$.
If $d(u)\geq 3$, $u\in V_{1} - v - v'$ and $u$ has a neighbor of degree $1$, then by Lemma~\ref{10l}, $d(u)\geq 3K$.  
By Lemma~\ref{13l}, $|N(u)-L(u)| - \frac{1}{3} |L(u)| \geq \frac{1}{3} d(u) \geq K  =15$. So, after giving away charges by (R2), $u$ keeps for itself charge at least $3$.
If $u \in \{v,v'\}$, then it originally had extra $\D+4$ of charge and it gives out by (R3) at most $\D + 4$.

If $u\in V_2$ and $d(u)\geq 4$, then by (R1), it keeps $4$ for itself.
Suppose finally that $u\in V_2$ and $d(u)=3$.  If it is adjacent to $v$ or $v'$, then by (R3), $ch^*(u)\geq 3 + 1 = 4$. 
Otherwise, 
by Remark~\ref{r1}, $u$ has a neighbor $y \in V_2$ with degree at least $K+1$ and by (R1) receives from $y$ charge $1-\frac{4}{K+1} > \frac{2}{3}$. \QED

%======================================================================
%Section:  Weak Vertices and Sponsors
%======================================================================
% 
%
\section{Weak Vertices and Sponsors}\label{sec:sponsors}
%
%
%======================================================================

A \emph{weak} vertex is either a $1$-vertex or a $2$-vertex with a neighbor of degree $2$. The \emph{sponsor}, $s(u)$, of a weak vertex $u$
is the unique neighbor of $u$ of degree at least $3$. By  Lemma~\ref{10l}, $d(s(u))\geq \frac{13}{5}K$ for each weak $u$.
A \emph{supersponsor} 
 is a vertex with at least
two neighbors that are weak. Notice that, for example, every donor is also a supersponsor.  By definition, each supersponsor is the sponsor for each of its weak neighbors.

\begin{lem}\label{15l}
Either $V_1$ or $V_2$ contains more than one supersponsor.
\end{lem}

\noindent
\textbf{Proof: } Suppose not. Choose $v_0 \in V_1$ and $w_0\in V_2$ so that no $x\in V(\mathbf{G})-v_0-w_0$ is a supersponsor.
For $x\in V(\mathbf{G})$, let $W(x)$ denote the set of weak neighbors of $x$. By our assumption, $|W(x)|\leq 1$ for each $x\in V(\mathbf{G})-v_0-w_0$.
Consider the following discharging. 

To start we let $ch(v_0)=d(v_0)+2\D+7$,  $ch(w_0)=d(w_0)+3$, $ch(u)=d(u)$ for each $u\in V(\mathbf{G})-v_0-w_0$. 
\begin{equation}\label{528}
\mbox{The total charge is $2(e_1+e_2+e_3+\D+5)$.}
\end{equation}
We redistribute charges 
according to the following rules.

(R1)  Each vertex $u$ of degree at least $4$ not adjacent to weak vertices gives to each neighbor
charge $\frac{d(u)-3}{d(u)}$ (and keeps $3$ for itself). 

(R2) Each vertex $u\in V(\mathbf{G})-v_0-w_0$ with $d(u)=3$ gives to each  neighbor of degree $2$ charge $1/4$.

(R3) Each sponsor $u \in V(\mathbf{G}) - v_{0} - w_{0}$ (then its degree is at least $\frac{13}{5}K$ by Lemma~\ref{10l}(b))
gives to each $x \in W(u)$ charge $2$ and to each other neighbor charge $\frac{d(u)-5}{d(u)}$, and leaves charge at least
$5-2\cdot|W(u)|\geq 3$ for itself. 

(R4)  Vertex $v_0$ gives $2$ to each neighbor and leaves $(2\D+d(v_0)+7)-2d(v_0)\geq 3$ 
for itself.

(R5)  Vertex $w_0$ gives $1$ to each neighbor and leaves 3
for itself.

We will show that the resulting charge, $ch^*(x)$, is at least $3$ for each $x\in V(\mathbf{G})$. Together with~\eqref{528}, this will contradict~\eqref{e1}.

Indeed, if $x$ is weak and has degree 1, then it must be in $V_1$ and so it will get $2$ by (R3) or by (R4). If it is weak and degree 2, then it gets at least $1$ by (R3), (R4), or (R5). If $d(x)=2$, and $x$ is not weak, then
$x$ gets at least $1-\frac{5\cdot 7}{13K}=1-\frac{7}{39}$ from its neighbor of degree at least $\frac{13K}{7}$ and at least $\frac{1}{4}$ from another neighbor; 
in total, more than $1$.
If $d(x)=3$, then $x$  gets at least $\frac{K-5}{K}=\frac{2}{3}$ from its neighbor of degree at least $K$,
 and gives away at most $\frac{2}{4}$ by (R2). Similarly, if $d(x)\geq 4$, then by (R1),(R3),(R4) or (R5), it reserves charge $3$ for itself.
\QED

\begin{lem}\label{16l}
If $V_i$ contains at least two supersponsors, then for each weak $w \in V_{3-i}$, the unique sponsor of $w$ is also contained in $V_{3-i}$.
\end{lem}

\noindent
\textbf{Proof: }
Suppose a weak   $w\in V_{3-i}$ is adjacent to a vertex $x_1\in V_i$ of degree at least $\frac{13}{5}K$. By Lemma~\ref{3l},
 $d(w)=2$ and $w$ has a neighbor $w'\in V_{3-i}$ with $d(w')=2$. Let $w''$ be the other neighbor of $w'$ (possibly, $w''\in V_i$).
 By the conditions of the lemma, there is a supersponsor $x_2 \in V_i - x_1$.  By Claim~\ref{12c}, there is a vertex $x_3\in V_i-N[x_2] - w''$ of degree at most $3$.  Send $x_2$ to $w$, $x_3$ to $w'$, and, if $w''\in V_{3-i}$, join $w''$ with the white neighbors of $x_3$ (there are at most $3$ of them) by yellow edges.  This way we eliminate all $d(x_2)+d(w)+1$ edges incident with $x_2$ or $w$ or $w'$, add at most $3$ yellow edges and increase $\D$ by at most $3$. Moreover, the remaining graph triple $\mathbf{G}'$ satisfies \eqref{e0} since for $i=1,2,3,$ \[\Delta_i\leq\Delta_i+3\leq (\D+4)+3\leq n+9-C<(n-2)-2.\] 
 Since $d(x_2) + d(w) + 1 \geq \frac{13}{5}K  + 3 \geq 18$, we see that $\partial (\mathbf{G}, \mathbf{G}') \geq 18-3-3 =12$. Hence, we are able to pack the remaining graph triple since $\mathbf{G}$ was a minimal counterexample.\QED

\begin{lem}\label{17l}
Each of $V_1$ and $V_2$ contains at least two supersponsors.
\end{lem}

\noindent
\textbf{Proof: } Suppose $V_{i}$ contains at most one supersponsor and this supersponsor is $w_0$, if exists.
Then by Lemma~\ref{15l}, $V_{3-i}$ contains two supersponsors $x_1$ and $x_2$. By Lemma~\ref{16l}, the sponsor of each weak vertex in $V_{i}$ is also in $V_{i}$. 
By Lemma~\ref{14l},  $\mathbf{G}$ has at most one donor.
Let $v_0$ denote such a vertex, if it exists. By~\eqref{dec243},  $v_0\in V_1$, and by definition it is a supersponsor.

\textbf{Case 1:} $i=2$. We use the following discharging. Let $ch(u) = d(u)$ for each $u \in V - v_{0} - w_{0}$. If $w_0$ and/or $v_0$ exist, then
let $ch(v_{0}) = d (v_{0}) + \Delta_{1} +  \Delta_{3|1} + 4$, and $ch(w_{0}) = d(w_{0}) +  \Delta_{2} +  \Delta_{3|2} + 4$.  By the definition of $\D$, 
\[\Delta_{1} +  \Delta_{3|1} + \Delta_{2} +  \Delta_{3|2}   \leq 2\D + 8, \] so the total charge is at most $2 (e_{1} + e_{2} + e_{3} + \D + 8)$.

Then we redistribute the charges using the following  set of rules.

(R1)  Each vertex $u$ of degree at least $5$ not adjacent to weak vertices gives to each neighbor
charge $\frac{d(u) - {19}/{6}}{d(u)} \geq \frac{1}{3}$ (and keeps $\frac{19}{6}$ for itself). 

(R2) Each vertex $u\in V(\mathbf{G})$ with $d(u) = 3$ or $d(u) = 4$ gives to each  neighbor of degree $2$ charge $\frac{1}{3}$.

(R3) Each sponsor $u\in V(\mathbf{G})$ (then  by Lemma~\ref{10l}(b) its degree is at least $\frac{13K}{5}=39$)
 but not a supersponsor gives charge 
$\frac{13}{6}$ to its  weak neighbor,  
and charges $\frac{d(u)-4.5}{d(u)}$ to each other neighbor.

(R4)  Each supersponsor $u\notin \{v_0,w_0\}$ gives $\frac{13}{6}$ to each adjacent $1$-vertex 
(by Lemma~\ref{14l} and the definition of $v_0$, there is at most $1$ such neighbor) and $\frac{d(u)-4.5}{d(u)}$ to each other neighbor.
 
 (R5) Each of $w_0$ and $v_0$ gives $\frac{11}{6}$ to each  neighbor.

 \medskip
We will show that the resulting charge, $ch^*(y)$, is at least $\frac{17}{6}$ for each $y\in V_1$ and
at least $\frac{19}{6}$ for each $y\in V_{2}$. This would mean  the total charge is at least $6n$, a contradiction to~\eqref{e1}.

Indeed, if $y$ is a 
$1$-vertex, then it is in $V_{1}$ and  will get $\frac{11}{6}$ by (R3), (R4),  or (R5).  
If $y$ is a
weak $2$-vertex and not adjacent to a supersponsor, then it will get $\frac{13}{6}$ from its sponsor  by (R3). 
If $y$ is a
weak $2$-vertex  adjacent to a supersponsor and $y \in V_1$, then by (R4) or (R5), it will get 
at least $1-\frac{4.5}{39} > \frac{5}{6}$ from its sponsor, and its resulting charge will be at least $\frac{17}{6}$ . 
If $y$ is a
weak $2$-vertex in $V_{2}$ adjacent to a supersponsor, then by Lemma~\ref{16l}, this supersponsor is $w_0$,
and $y$ gets $\frac{11}{6}$ from $w_0$.

If $d(y)=2$, and $y$ is not weak, then by Lemma~\ref{10l}(a), $y$ has a neighbor of degree at least $\left\lceil\frac{13K}{7}\right\rceil=28$.
So
$y$ gets from it at least $1-\frac{4.5}{28}$  (by (R1), (R3), (R4) or (R5))
 and at least $\frac{1}{3}$ from another neighbor (by one of (R1)--(R5)). Then $ch^*(y)\geq 3-\frac{4.5}{28} + \frac{1}{3}> \frac{19}{6}$. 
If $d(y) = 3$ and $y$ has two neighbors of degree $2$, then 
by Lemma~\ref{10l}(b), $y$ has a neighbor $x$  of degree at least $\frac{13K}{5}=39$, so it gets from $x$
  at least $\frac{39-4.5}{39}\geq \frac{5}{6}$,
 and gives away at most $\frac{2}{3}$ by (R2). If $d(y) = 3$ and $y$ has at most one neighbor of degree $2$, then it gets from its neighbor
 of degree at least $\left\lceil\frac{13K}{10}\right\rceil=20$ charge at least $\frac{15.5}{20}$ and gives away at most $\frac{1}{3}$.
 If $d(y) = 4$, then $y$ gets at least $\frac{K-5}{K}=\frac{2}{3}$ from it neighbor of degree at least $K$ and gives away at most $3 \cdot \frac{1}{3} = 1$
 by (R2).
 If $d(y) \geq 5$ and $y$ has no weak neighbors, then it leaves $\frac{19}{6}$ for itself by (R1).
 
 If $y$ has a weak neighbor and  $y \notin\{v_0,w_0\}$, then $d(y)\geq 39$ and
  by  (R3) or (R4), it reserves for itself charge 
  \[d(y)-\frac{13}{6}-(d(y)-1)\frac{d(y)-4.5}{d(y)}=-\frac{13}{6}+\frac{5.5d(y)-4.5}{d(y)}=\frac{10}{3}-\frac{4.5}{d(y)}\geq \frac{10}{3}-\frac{4.5}{39}>\frac{19}{6}. \]
 The vertex $w_{0}$ gives away charge $\frac{11}{6} d_{2}(w_{0}) + \frac{11}{6} d_{3}(w_{0}) \leq d(w_{0}) +  \Delta_{2} + \Delta_{3|2}$ and saves 
 more than $4$ for itself.  Similarly, $v_{0}$ saves 
 more than $4$ for itself. This proves the case.

\textbf{Case 2:} $i=1$. In this case either $v_0$ does not exist, or $v_0=w_0$. The discharging is very similar to that
in Case 1, but a bit simpler. 
 Let $ch(u) = d(u)$ for each $u \in V - w_{0}$. If $w_0$ exists, then
let $ch(w_{0}) = d (w_{0}) + 2\D+ 4$. So, the total charge is at most $2 (e_{1} + e_{2} + e_{3} + \D + 4)$.
The first $3$ rules of discharging are again (R1)--(R3), but instead of (R4) and (R5), we have

(Q4)  Each supersponsor $u\neq w_0$ gives  $\frac{d(u)-4.5}{d(u)}$ to each  neighbor.
 
 (Q5) Vertex $w_0$  gives $\frac{13}{6}$ to each  neighbor.

 \medskip
Symmetrically to Case 1, we will show that the resulting charge, $ch^*(y)$, is at least $\frac{19}{6}$ for each $y\in V_1$ and
at least $\frac{17}{6}$ for each $y\in V_{2}$, again yielding a contradiction to~\eqref{e1}.
 
 If $y$ is a 
$1$-vertex, then it is in $V_{1}$ and its neighbor also is in $V_1$. Since all supersponsors apart from $w_0$ are in $V_2$,
Rule (Q4) does not apply to $y$, so $y$ will get $\frac{13}{6}$ by (R3)  or (Q5).  
If $y$ is a
weak $2$-vertex and not adjacent to a supersponsor, then it will get $\frac{13}{6}$ from its sponsor  by (R3). 
If $y$ is a
weak $2$-vertex  adjacent to a supersponsor and $y \in V_2$, then by (Q4) or (Q5), it will get 
at least $1-\frac{4.5}{13K/5} = 1-\frac{3}{26}$ from its sponsor, so that its resulting charge will be more than $\frac{17}{6}$ . 
If $y$ is a
weak $2$-vertex in $V_{1}$ adjacent to a supersponsor, then by Lemma~\ref{16l}, this supersponsor is $w_0$,
and $y$ gets $\frac{13}{6}$ from $w_0$.
 
 Counting of charges for other vertices apart from $w_0$ simply repeats that in Case 1 (using (Q4) and (Q5) in place of (R4) and (R5)).
 Since the starting charge of $w_0$ was at least $3d(w_0)$, by (Q5), its new charge is at least $\frac{5}{6}d(w_0)+4>4$.
\QED

%======================================================================
%Section:  End of the Proof
%======================================================================
% 
%
\section{Proof of Theorem~\ref{thm:detailed main}}\label{sec:end}
%
%
%======================================================================
%====================================================================== 

By Lemma~\ref{17l},  $V_{1}$ contains  supersponsors $x_1$ and $x_2$ and $V_{2}$ contains  supersponsors $y_1$ and $y_2$. Let $v_1$ (resp. $w_1$) be a weak neighbor of $x_1$ (of $y_1$), let $v'_1$ ($w'_1$) be the other neighbor of it  which is of degree $2$ if it exists, and let $v''_1$  ($w''_1$) be the other neighbor of $v'_1$  (of $w'_1$). Let $v_2$ ($w_2$) be a weak neighbor of $x_2$ (of $y_2$) that is \emph{not adjacent} to $v_1$ (to $w_1$); this is possible since $x_{2}$ ($y_{2}$) is adjacent to multiple weak vertices.  Let $v'_2$ ($w'_2$) be the other neighbor of it which is again of degree $2$ if it exists, and let $v''_2$  ($w''_2$) be the other neighbor of $v'_2$  (of $w'_2$).
  
We are now ready to construct our packing.  For $j = 1,2$, begin by placing $x_{j}$ on $w_{j}$, and $v_{j}$ on $y_{3-j}$.  Notice that  by  Lemma~\ref{16l}, $v_{j} \in V_{1}$ and $w_{j} \in V_{2}$ so this assignment is well defined.  Since the weak vertices have only one sponsor, $v_{j}$ is not adjacent to $x_{3-j}$, $y_{1}$, nor $y_{2}$, and $w_{j}$ is not adjacent to $y_{3-j}$, $x_{1}$, nor $x_{2}$.  Together with the fact that $v_{1}$ ($w_{1}$) was chosen to be not adjacent to $v_{2}$ ($w_{2}$), we see that these mappings do not violate the packing property.
 
As we extend this packing, we only need to ensure that $v_{j}'$ is not mapped to a vertex in $N_{2}(y_{3-j})$ and no vertex in $N_{1}(x_{j})$ is mapped to $w_{j}'$.  This can only be an issue if  $v_{j}' \in V_{1}$ ($w_{j}' \in V_{2}$) and in this case, we will find an appropriate assignment for $v_{j}'$.  If $v_{j}' \in V_{2}$ ($w_{j}' \in V_{1}$), we will simply ignore this part of the construction.

By Claim~\ref{12c}, there is a vertex $x'_1\in V_1-N(x_1)-\bigcup_{i=1,2}\{v_i,v_i',v_i'',w_i,w_i',w_i''\}$
 ($y'_1\in V_2-N(y_1)-\bigcup_{i=1,2}\{v_i,v_i',v_i'',w_i,w_i',w_i''\}$)
 with degree at most $3$.  Similarly,
 there are vertices $x'_2\in V_1-N(x_2)-x'_1-\bigcup_{i=1,2}\{v_i,v_i',v_i'',w_i,w_i',w_i''\}$
 and  $y'_2\in V_2-N(y_2)-y'_1 - \bigcup_{i=1,2}\{v_i,v_i',v_i'',w_i,w_i',w_i''\}$
  of degree at most $3$.

\begin{figure}[h!]
\begin{center}
\includegraphics[width=.4\textwidth]{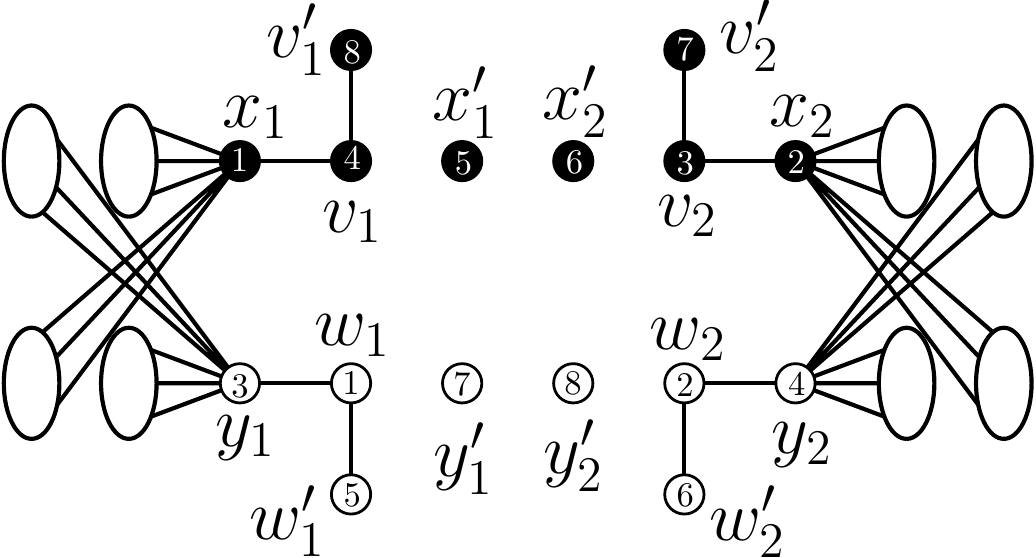}
\end{center}
\caption{Sketch of Packing}
\label{fig:EndProof}%
\end{figure}

For the following mappings, refer to Figure \ref{fig:EndProof}. If $w_{j}' \in V_{2}$, then send $x'_j$ to $w_j'$ and, if $w''_j\in V_2$, add the yellow edges connecting $w''_j$ with the at most $3$ white neighbors of $x'_j$.  Similarly, if $v_{j}' \in V_{1}$, then send $v_{j}'$ to $y'_{3-j}$ (if $v_j'\in V_1$) and, if $v''_j\in V_1$, add the yellow edges connecting $v_{j}''$ with the at most three white neighbors of $y_{3-j}'$.

Let $\mathbf{G}'$ be the triple obtained by deleting the assigned vertices. By construction, if $\mathbf{G}'$ packs, then together with our placement, we get a packing of $\mathbf{G}$.  We decreased $n$ by at most $8$ and decreased the number of edges by at least $d(x_1)+d(x_2)+d(y_1)+d(y_2) - 16 \geq 12K - 16$. We have increased $\D$ by at most $6$ (with the new yellow edges). So, $\partial (\mathbf{G},\mathbf{G}') \geq 12K - 22 \geq 24 = 3 (n - n')$.  Since $d_{i} (v) \leq \D + 4 \leq n - C +6$ for every $v \in V$ (and $C \geq 8$),~\eqref{e0} holds for $\mathbf{G}'$. Thus $\mathbf{G}'$ (and hence $\mathbf{G}$) packs, a contradiction to the choice of $\mathbf{G}$.\QED

\textbf{Case 1:}  \emph{The vertices $w_{0} \in V_{2}$ and $v_{0} \in V_{1}$ are distinct.}  In this case, $w_{0} \in V_{2}$ is the only supersponsor in $V_{2}$.  

\textbf{Case 2:} \emph{The vertex $v_{0}$ does not exist or $w_{0} = v_{0}$.}  In this case,  the initial charge will be slightly different.  
For each $u \in V - w_{0}$, $ch(u) = d(u)$ and $ch(w_{0} ) = d(w_{0}) + 2\D + 16$.  As in Case 1,  the total charge is at most $2(e_{1} + e_{2} + e_{3} + \D + 8)$.  Further, the charge assigned to $w_{0}$ in this case is at least the charge assigned to it in Case 1.

\textbf{Case 3:} \emph{The vertex $v_{0}$ exists but $w_{0}$ does not.}  This case is symmetric to Case 2.  For each $u \in V - v_{0}$, $ch(u) = d(u)$ and $ch(v_{0} ) = d(v_{0}) + 2\D + 16$.  As in the previous cases, the total charge is at most $2(e_{1} + e_{2} + e_{3} + \D + 8)$.  Further, the charge assigned to $v_{0}$ is at least the charge assigned to it in Case 1.

For \emph{all} cases, we redistribute the charges using the following same set of rules.

(R1)  Each vertex $u$ of degree at least $5$ not adjacent to weak vertices gives to each neighbor
charge $\frac{d(u) - \frac{19}{6}}{d(u)} \geq \frac{1}{3}$ (and keeps $\frac{19}{6}$ for itself). 

(R2) Each vertex $u\in V(\mathbf{G})$ with $d(u) = 3$ or $d(u) = 4$ gives to each  neighbor of degree $2$ charge $\frac{1}{3}$.

(R3) Each non-weak vertex $u\in V(\mathbf{G})$ adjacent to a weak vertex (then its degree is at least $K$ by Lemma~\ref{10l}) but not a supersponsor gives charge 
$\frac{11}{6}$ to its neighbor of degree $1$ (if such neighbor exists) or $\frac{7}{6}$ to its weak neighbor of degree $2$,  
and charges $\frac{d(u)-5}{d(u)}$ to each other neighbor.

(R4)  Each supersponsor $u\notin \{v_0,w_0\}$ gives $\frac{11}{6}$ to each adjacent $1$-vertex 
(by Lemma~\ref{14l} and the definition of $v_0$, there is at most $1$ such neighbor) and $\frac{d(u)-5}{d(u)}$ to each other neighbor.
 
 (R5) The vertex $w_0$ gives $\frac{11}{6}$ to each  neighbor.
 
 (R6) The vertex $v_{0}$, if it is distinct from $w_{0}$, gives charge $\frac{11}{6}$ to each  neighbor.

\begin{rem}\label{r2}
If in the statement of Lemma~\ref{10l}, $v \in V_{i} $, $2\leq d(v)=t \leq 4$ and at least one neighbor of $v$ has degree less than $5$, then 
either $v$ has a neighbor in $V_i$ of degree at least $\frac{13K}{3t-1}$, or $v$ is adjacent to all vertices in $V_{3-i}$ of degree at least~$F$.
\end{rem}

\noindent
\textbf{Proof: }  Let $N_1(v):=
 \{v_1,\ldots,v_s\}$. If $s<t$, then the proof of Lemma~\ref{10l} works. So suppose $s=t$ and $d(v_s)\leq 4$.
We almost word by word repeat the proof of Lemma~\ref{10l} with $\frac{13K}{3t-1}$ in place of $\frac{13K}{3t+1}$, only the number of added yellow edges is now at most
$3 \left((s-1) (\frac{13K}{3t-1}-1)+4\right)$, so that instead of~\eqref{dec29}, we get
\[\partial (\mathbf{G},\mathbf{G}') \geq d(v) + d(w)  - (3s-2) (\frac{13K}{3t-1}-1)-12 \]
\[\geq (s+1) + d(w) -  13(1-\frac{1}{3t-1})K+(3s-2)-12.\]
Since $d(w)\geq F=13K$ and $2\leq s=t\leq 4$, this is at least 
\[  F - 13 K+\frac{13K}{11}+(4s-1)-12 > 3s+3.\]
Thus, as in the proof of Lemma~\ref{10l}, $\mathbf{G}'$ packs and so $\mathbf{G}$ packs.\QED

\bibliographystyle{abbrv}
\bibliography{ListPacking}
\end{document}